\documentclass[12pt]{article}
\usepackage{amsthm,amsfonts,amsmath}
\newtheorem{theorem}{Theorem}[section]

\newtheorem{lemma}[theorem]{Lemma}
\newtheorem{proposition}[theorem]{Proposition}
\newtheorem{conjecture}[theorem]{Conjecture}

\theoremstyle{definition}

\theoremstyle{remark}
\newtheorem{remark}{Remark}[section]

\newcommand{\1}{{{\mathchoice {\rm 1\mskip-4mu l} {\rm 1\mskip-4mu l}
{\rm 1\mskip-4.5mu l} {\rm 1\mskip-5mu l}}}}
\newcommand{\dslash}{/\mskip-6mu/}

\newcommand{\Z}{{\mathbb{Z}}}
\newcommand{\R}{{\mathbb{R}}}
\newcommand{\C}{{\mathbb{C}}}

\newcommand{\Aa}{{\mathcal{A}}}   
\newcommand{\Bb}{{\mathcal{B}}}
\newcommand{\Dd}{{\mathcal{D}}}
\newcommand{\Ee}{{\mathcal{E}}}

\newcommand{\Gg}{{\mathcal{G}}}   

\newcommand{\Jj}{{\mathcal{J}}}

\newcommand{\Ll}{{\mathcal{L}}}   
\newcommand{\Mm}{{\mathcal{M}}}   

\newcommand{\Oo}{{\mathcal{O}}}

\newcommand{\Ss}{{\mathcal{S}}}

\newcommand{\Ww}{{\mathcal{W}}}
\newcommand{\Xx}{{\mathcal{X}}}

\newcommand{\im}{{\rm im }}        
\newcommand{\tr}{{\rm tr }}        
\newcommand{\id}{{\rm id}}         

\newcommand{\INDEX}{{\rm index}}   
 
\newcommand{\grad}{{\rm grad }}    
\newcommand{\Lie}{{\rm Lie}}          
\newcommand{\Vect}{{\rm Vect}}        
\newcommand{\Vol}{{\rm Vol}}          
\newcommand{\Ham}{{\rm Ham}}          
\newcommand{\Per}{{\rm Per}}          
\newcommand{\Map}{{\rm Map}}          
\newcommand{\End}{{\rm End}}          
\newcommand{\HF}{{\rm HF}}            

\newcommand{\ad}{{\rm ad}}

\newcommand{\point}{{\rm pt}}
\newcommand{\ev}{{\rm ev}}

\newcommand{\FLAT}{{\rm flat}}
\newcommand{\dvol}{{\rm dvol}}
\newcommand{\loc}{{\rm loc}}

\newcommand{\eps}{{\varepsilon}}

\renewcommand{\i}{{\iota}}

\newcommand{\om}{{\omega}}

\newcommand{\Om}{{\Omega}}
\newcommand{\G}{{\rm G}}
\newcommand{\EG}{{\rm EG}}
\newcommand{\BG}{{\rm BG}}

\newcommand{\SO}{{\rm SO}}
\newcommand{\U}{{\rm U}}
\newcommand{\SU}{{\rm SU}}

\newcommand{\PSL}{{\rm PSL}}

\newcommand{\so}{{\mathfrak s \mathfrak o}}
\renewcommand{\u}{{\mathfrak u}}

\newcommand{\g}{{\mathfrak g}}         
\newcommand{\Cinf}{C^{\infty}}

\newcommand{\SW}{{\rm SW}}
\newcommand{\Gr}{{\rm Gr}}

\newcommand{\Sreg}{{\mathcal{S}}_{\rm reg}}

\newcommand{\inner}[2]{\langle #1, #2\rangle}

\newcommand{\winner}[2]{\langle #1{\wedge}#2\rangle}   
\def\NABLA#1{{\mathop{\nabla\kern-.5ex\lower1ex\hbox{$#1$}}}}
\def\Nabla#1{\nabla\kern-.5ex{}_{#1}}
\def\Tabla#1{\Tilde\nabla\kern-.5ex{}_{#1}}
\renewcommand{\Tilde}{\widetilde}
\renewcommand{\Hat}{\widehat}

\newcommand{\p}{{\partial}}

\newcommand{\IMP}{\Longrightarrow}

\newcommand{\INTO}{\hookrightarrow}
\newcommand{\TO}{\longrightarrow}
\title{$J$-holomorphic curves, moment maps, \\
       and invariants of Hamiltonian group actions}

\author{Kai~Cieliebak\\
        Stanford University
        \and
        Ana Rita Gaio \\
        University of Warwick
        \and
        Dietmar~A.~Salamon\\ 
        ETH-Z\"urich}

\date{17 September 1999}

\begin{document}

\maketitle



\tableofcontents



\section{Introduction} 

This paper outlines the construction 
of invariants of Hamiltonian group 
actions on symplectic manifolds.
The invariants are derived from the solutions 
of a nonlinear first order elliptic partial differential
equation involving the Cauchy-Riemann operator,
the curvature, and the moment map (see~(\ref{eq:new}) below).
They are related to the Gromov invariants 
of the reduced spaces.

Our motivation arises from the proof of the Atiyah-Floer
conjecture in~\cite{DOSAL1,DOSAL2,DOSAL3} which deals with
the relation between holomorphic curves $\Sigma\to\Mm_S$
in the moduli space $\Mm_S$ of flat connections over 
a Riemann surface $S$ and anti-self-dual instantons
over the $4$-manifold $\Sigma\times S$.  
In~\cite{AB} Atiyah and Bott interpret the space 
$\Mm_S$ as a symplectic quotient of the space 
$\Aa_S$ of connections on $S$ by the action 
of the group $\Gg_S$ of gauge transformations.  
A moment's thought shows that the various terms in the 
anti-self-duality equations over $\Sigma\times S$
(see equation~(\ref{eq:asd}) below) can be interpreted 
symplectically.   Hence they should give rise to meaningful
equations in a context where the space $\Aa_S$ is replaced
by a finite dimensional symplectic manifold $M$ and the gauge
group $\Gg_S$ by a compact Lie group $\G$ with
a Hamiltonian action on $M$.  In this paper we show
how the resulting equations give rise to invariants 
of Hamiltonian group actions.  The same adiabatic 
limit argument as in~\cite{DOSAL3} then leads
to a correspondence between these invariants
and the Gromov--Witten invariants of the 
quotient $M\dslash\G$ (Conjecture~\ref{con:adiabatic}).  
This correspondence is the subject of the 
PhD thesis~\cite{GAIO} of the second author.  

In Section~\ref{sec:background} we review the relevant 
background material about Hamiltonian group actions,
gauge theory, equivariant cohomology, and holomorphic
curves in symplectic quotients. The heart of this paper 
is Section~\ref{sec:new}, where we discuss the equations
and their properties, outline the construction 
of the invariants, and indicate several potential
applications.  One interesting point is a result
about the compactness of the moduli spaces
(Proposition~\ref{prop:compact}) which has
no analogue for moduli spaces of holomorphic curves.
Hence the invariants should lead in many cases to
a definition of the Gromov--Witten invariants 
over the integers.  
Via the adiabatic limits and wallcrossing arguments,
the invariants should also give rise
to relations between the Gromov--Witten invariants
of symplectic quotients at different values
of the moment map.  This is reminiscent of the 
work of Martin~\cite{MARTIN1,MARTIN2,MARTIN3}
about the ordinary cohomology of symplectic quotients. 
Section~\ref{sec:floer} deals with the corresponding Floer theory
and Section~\ref{sec:ex} discusses several examples. 

A few historical comments are in order.  The last author 
discovered the equations~(\ref{eq:new}) and some of their 
potential applications in 1995 and gave lectures about them
at various places.  But until now this work did not appear
in preprint form. The second author has been working 
on her thesis on this subject since 1996.
Recently, Ignasi Mundet discovered 
the same equations independently, starting from a different 
angle, and in his thesis~\cite{MUNDET} developed a programme 
along similar lines as outlined in the present paper. 
For $M=\C^n$ the equations also appeared in the physics 
literature (starting from Witten's work in~\cite{WITTEN})
where they are known as {\it gauged sigma models}. 


\section{Background}\label{sec:background}


\subsection{Hamiltonian group actions}\label{sec:moment}

Let $(M,\om)$ be a symplectic manifold
and $\G$ be a compact connected Lie group
that acts on $M$ by symplectomorphisms.
Let $\g=\Lie(\G)$ denote the Lie algebra of $\G$
and, for every $\xi\in\g$, denote by $X_\xi:M\to TM$
the vector field whose flow is given by the action 
of the one-parameter subgroup generated by $\xi$. 
We assume throughout that the Lie algebra $\g$ 
carries an invariant inner product $\inner{\cdot}{\cdot}$.
The action of $\G$ on $M$ is called {\bf Hamiltonian}
if there exists an equivariant function 
$\mu:M\to\g$ such that, for every $\xi\in\g$,
\begin{equation}\label{eq:moment}
     d\inner{\mu}{\xi} = \i(X_\xi)\om.
\end{equation}
This means that $X_\xi$ is the Hamiltonian vector field
of the function $\inner{\mu}{\xi}$.  
The function $\mu$ is called a {\bf moment map}. 

Suppose that $\tau$ is a regular value of $\mu$
and that the isotropy subgroup 
$$
    \G_\tau := \left\{g\in\G\,|\,\tau=g\tau g^{-1}\right\}
$$
acts freely on $\mu^{-1}(\tau)$.  
Then the {\it Marsden-Weinstein quotient}
$$
     M\dslash\G(\tau) := \mu^{-1}(\tau)/\G_\tau
$$
is a smooth manifold and it inherits the symplectic 
structure from $M$. To be more precise, 
for $x\in\mu^{-1}(\tau)$,
there is a chain complex
\begin{equation}\label{eq:exact-mu}
    0\TO\g_\tau
     \stackrel{L_x}{\TO} T_xM
     \stackrel{d\mu(x)}{\TO} \g
     \TO 0,
\end{equation}
where $L_x:\g\to T_xM$ is defined by $L_x\xi=X_\xi(x)$
and 
$$
    \g_\tau=\{\xi\in\g\,|\,[\xi,\tau]=0\}
$$
is the Lie algebra of $\G_\tau$.
The image of $L_x$ is the symplectic complement
of the kernel of $d\mu(x)$. 
Moreover, the formula
$$
    d\mu(x)L_x\xi = [\xi,\mu(x)]
$$
shows that $\im L_x\cap\ker d\mu(x)=L_x\g_\tau$. 
Thus the quotient $\ker d\mu(x)/L_x\g_\tau$ inherits 
the symplectic structure of $T_xM$, and it can be identified
with the tangent space of $M\dslash\G(\tau)$ at $[x]$. 

\begin{remark}\label{rmk:quotient}\rm
{\bf (i)}
An orbit $\Oo\subset\g$ under the adjoint action 
admits a natural K\"ahler structure.  
The tangent space of $\Oo$ at $\tau$ is
$$
     T_\tau\Oo = \{[\xi,\tau]\,|\,\xi\in\g\} = \g_\tau^\perp
$$
and the symplectic form is given by
$$
     \sigma_\tau([\xi,\tau],[\eta,\tau]) := \inner{\tau}{[\xi,\eta]}.
$$
An explicit formula for the complex structure uses 
the decomposition of $T_\tau\Oo$ 
into the eigenspaces
$$
     V_\alpha 
     = \left\{\xi\in\g\,|\,
       [\tau,[\tau,\xi]]=-\alpha^2\xi\right\}
$$
for $\alpha>0$.  Let $\g\to V_\alpha:\xi\mapsto\xi_\alpha$
denote the orthogonal projection onto $V_\alpha$ and define
$A_\tau:\g\to\g$ by 
$$
     A_\tau\xi = \sum_\alpha \alpha\xi_\alpha
$$
where the sum runs over all $\alpha>0$ such that 
$-\alpha^2$ is an eigenvalue of $\ad(\tau)^2$.     
Then the complex structure on $T_\tau\Oo$ is given by
$$
     J_\tau\xi = \sum_\alpha\frac{1}{\alpha}[\tau,\xi_\alpha],\qquad
     J_\tau[\xi,\tau] = A_\tau\xi.
$$ 
The adjoint action of $\G$ on $\Oo$ is Hamiltonian 
and one checks easily that the inclusion 
$\Oo\INTO\g$ is a moment map for this action.

\smallskip
\noindent{\bf (ii)}
If $\tau\in\Oo$ then the quotient 
$M\dslash\G(\tau)$ can be naturally identified with the 
symplectic quotient of $M\times\Oo$ at the zero value of 
the moment map.  Here the product $M\times\Oo$ 
is equipped with the symplectic
form $\om-\sigma$ and the moment map 
$
     \mu_\Oo:M\times\Oo\to\g
$
is given by
$
     \mu_\Oo(x,\tau)=\mu(x)-\tau.
$
Hence 
$
    (M\times\Oo)\dslash\G
    \cong\mu^{-1}(\Oo)/\G
    \cong\mu^{-1}(\tau)/\G_\tau
    = M\dslash\G(\tau)
$
for $\tau\in\Oo$. If $\mu(x)=\tau$ then
$$
    T_{[x,\tau]}(M\times\Oo)\dslash\G
    = \frac
    {\left\{(v,\eta)\in T_xM\times T_\tau\Oo\,|\,
    d\mu(x)v=\eta\right\}}
    {\left\{(L_x\xi,[\xi,\tau])\,|\,\xi\in\g\right\}},
$$
Equivalently, this tangent space can be identified with the 
space of all pairs $(v,\eta)\in T_xM\times T_\tau\Oo$
that satisfy
\begin{equation}\label{eq:TMO}
    d\mu(x)v=\eta,\qquad {L_x}^*v+J_\tau\eta = 0,
\end{equation}
where the adjoint operator ${L_x}^*$ is understood
with respect to a $\G$-invariant inner product that arises
from a $\G$-invariant almost complex structure $J$ on $M$
that is compatible with $\om$.  
On the other hand,
$$
    T_{[x]}M\dslash\G(\tau)
    = \frac{\ker d\mu(x)}{\left\{L_x\xi\,|\,[\xi,\tau]=0\right\}}.
$$
The harmonic representative of a tangent vector 
$v\in\ker d\mu(x)$ is given by 
\begin{equation}\label{eq:harmonic}
    \pi(v) = (v+L_x\xi,[\xi,\tau]),\qquad
    {L_x}^*L_x\xi + A_\tau\xi + {L_x}^*v=0.
\end{equation}
The formula
$
    \om(v,w)
    = \om(v+L_x\xi,w+L_x\eta)
      - \inner{\tau}{[\xi,\eta]}
$
for $v,w\in\ker\,d\mu(x)$ shows that the
symplectic forms agree. 

\smallskip
\noindent{\bf (iii)}
Assume that $(M,\om,J)$ is a K\"ahler manifold and that the action
of $\G$ preserves the K\"ahler structure. Then the 
$\G$-action extends to an action of the complexified
group $\G^c$ that preserves the complex structure (cf.~\cite{GS1}).
The extended action does not preserve the K\"ahler form.   
Suppose that $\tau\in\g$ is a central element and
denote 
$$
     M^\tau
     = \left\{x\in M\,|\,\exists g\in\G^c
           \mbox{ such that }\mu(gx)=\tau\right\}.
$$
Then the complex quotient $M^\tau/\G^c$ can be naturally identified 
with $M\dslash\G(\tau)$.  This means that if 
$\mu(gx)=\mu(x)=\tau$ for $g\in\G^c$ and $x\in M$, 
then $gx$ lies in the $\G$-orbit of $x$. 
The proof of this observation relies on the fact that 
any $g\in\G^c$ can be written in the form
$g=\exp(i\eta)h$ where $h\in\G$ and $\eta\in\g$. 
Now consider the path $[0,1]\to M:x(t)=\exp(it\eta)hx$
running from $x(0)=hx$ to $x(1)=gx$.   This path satisfies
$$
     \dot x(t) = JX_\eta(x(t))
$$
and hence 
$$
     \frac{d}{dt}\inner{\mu(x(t))}{\eta}
         = \inner{d\mu(x(t))JX_\eta(x(t))}{\eta}
         = \left|X_\eta(x(t))\right|^2.
$$
The last identity follows from the definition
of the moment map. Since $\tau$ is a central element
we have $\mu(x(0)) = \mu(x(1))$, hence $x(t)$ is independent of $t$,  
and hence $gx=hx\in\G x$.

\smallskip
\noindent{\bf (iv)}
The study of complex quotients of the form $M^\tau/\G^c$
and their relation to the Marsden-Weinstein quotients
is the subject of {\it geometric invariant theory}
(cf.~\cite{MFK,DONALDSON1,DONALDSON5,DONALDSON6}). 
In his beautiful recent paper~\cite{DONALDSON6} Donaldson
treats the infinite dimensional case where $\G$ 
is replaced by the group of volume preserving 
diffeomeorphisms of a manifold $S$, and $M$ by the
space of maps from $S$ to $M$.  
In~\cite{DONALDSON5} Donaldson 
discusses another interesting case where 
$\G$ is the group of symplectomorphisms and $M$ is the 
manifold of almost complex structures compatible
with the given symplectic form. 

\smallskip
\noindent{\bf (v)}
The condition that $\G_\tau$ acts freely
on $\mu^{-1}(\tau)$ is rather strong.  
In general, if $\tau$ is a regular value of 
$\mu$, then the action of $\G_\tau$ 
on $\mu^{-1}(\tau)$ has finite isotropy subgroups
and the quotient $M\dslash\G(\tau)$ 
is a symplectic orbifold. 
Much of the dicussion in this
paper extends to that case. 
\end{remark}


\subsection{Connections and curvature}\label{sec:gauge}

Let $X$ be a compact oriented smooth manifold
and $P\to X$ be a principal $\G$-bundle. 
We think of $\G$ as acting on $P$ on the right
and denote the infinitesimal action 
by $p\xi\in T_pP$ for $p\in P$ and $\xi\in\g$.
A {\bf connection} on $P$ is an equivariant
horizontal subbundle of $TP$.  Any such 
subbundle determines an equivariant
$1$-form $A\in\Om^1(P,\g)$
whose kernels are the horizontal subspaces
and which identifies the vertical subspaces with $\g$. 
Thus 
$$
     A_{ph}(vh) = h^{-1}A_p(v)h,\qquad
     A_p(p\xi) = \xi
$$
for $p\in P$, $v\in T_pP$, $h\in\G$, and $\xi\in\g$.  
A $1$-form $A\in\Om^1(P,\g)$ that satisfies these conditions 
is called a {\bf connection $1$-form} and 
the space of connection $1$-forms will be denoted by 
$\Aa(P)$.  A {\bf gauge transformation} of $P$ is a 
smooth function $g:P\to\G$ that is equivariant with respect 
the adjoint action of $\G$ on itself, i.e.
$
     g(ph) = h^{-1}g(p)h
$
for $p\in P$ and $h\in\G$. The group of gauge transformations
will be denoted by $\Gg=\Gg(P)$.  It acts on $\Aa(P)$ 
contravariantly by 
$$
     g^*A = g^{-1}dg + g^{-1}Ag.
$$
Thus $g^*A$ is the pullback of $A$ under the automorphism
$P\to P:p\mapsto pg(p)$.   Let $\Om^k_\ad(P,\g)$ denote the 
space of equivariant and horizontal $k$-forms on $P$
with values in $\g$.  Any such form descends to a $k$-form 
on $X$ with values in the adjoint bundle
$\g_P:=P\times_\ad\g$.   Every connection 
$A\in\Aa(P)$ gives rise to a covariant derivative operator
$d_A:\Om^k_\ad(P,\g)\to\Om^{k+1}_\ad(P,\g)$ given by 
$$
     d_A\alpha =  d\alpha + [A\wedge\alpha].
$$
It is interesting to note that 
$\Om^0_\ad(P,\g)$ is the Lie algebra of the gauge group, 
$\Om^1_\ad(P,\g)$ is the tangent space of the space of connections, 
and the infinitesimal action of $\Lie(\Gg(P))$ 
on $\Aa(P)$ is given by the covariant derivative.  

Now suppose that $X=\Sigma$ is a compact Riemann surface.
Then the space $\Aa(P)$ carries
a natural symplectic form 
$$
    \Om(\alpha,\beta) = \int_\Sigma\winner{\alpha}{\beta}.
$$
Atiyah and Bott~\cite{AB} noted that the action of $\Gg(P)$
on $\Aa(P)$ is Hamiltonian and that a moment map
is given by the curvature 
$$
    F_A = dA + \frac12[A\wedge A]\in \Om^2_\ad(P,\g).
$$
Thus the Marsden-Weinstein quotient is the moduli
space of gauge equivalence classes of flat connections on $P$. 
The analogue of the chain complex~(\ref{eq:exact-mu}) 
for $\tau=0$ in gauge theory is given by
\begin{equation}\label{eq:exact-A}
    0\TO\Om^0_\ad(P,\g)
     \stackrel{d_A}{\TO} \Om^1_\ad(P,\g)
     \stackrel{d_A}{\TO} \Om^2_\ad(P,\g)
     \TO 0.
\end{equation}
Here $d_A:\Om^0_\ad\to\Om^1_\ad$ is the infinitesimal
action of the gauge group and 
$d_A:\Om^1_\ad\to\Om^2_\ad$ is the differential of the function
$\Aa(P)\to\Om^2_\ad:A\mapsto F_A$. 
The formula 
$
    d_Ad_A\tau=[F_A\wedge\tau]
$
shows that $d_A\circ d_A=0$ whenever $A$ is flat.


\subsection{Equivariant cohomology}\label{sec:eco}

Let $\EG$ be a contractible space on which the group 
$\G$ acts freely.  Unless $\G=\{\1\}$ this space is
necessarily infinite dimensional.
The {\bf equivariant (co)homology} 
of a $\G$-space $M$ is defined by 
$$
    H_\G^*(M;R) = H^*(M\times_\G\EG;R),\qquad
    H^\G_*(M;R) = H_*(M\times_\G\EG;R).
$$
Since there is a natural projection 
$
    M\times_\G\EG\to \EG/\G=:\BG,
$
$H^*_\G(M;R)$ is a module over $H^*(\BG;R)$.
Explicit representatives of equivariant homology
classes can be constructed as follows. 
Let $X$ be a compact oriented smooth 
$k$-manifold (without boundary)
and $\pi:P\to X$ be a principal $\G$-bundle.  
An equivariant map $u:P\to M$ 
determines an equivariant homology class 
$$
    [u]=u^\G_*(\pi^\G_*)^{-1}[X]\in H^\G_k(M;\Z).
$$
Here $u^\G_*:H^\G_*(P;\Z)\to H^\G_*(M;\Z)$ and
$\pi^\G_*:H^\G_*(P;\Z)\to H_*(X;\Z)$ denote 
the induced homomorphisms on equivariant homology
and $[X]\in H_k(X;\Z)$ denotes the fundamental class.
Since $\G$ acts freely on $P$, $\pi^\G_*$ is an isomorphism.

\begin{remark}\label{rmk:P}\rm
For every principal bundle $P\to X$ 
there exists an equivariant map $\phi:P\to\EG$.  
The map $P\to M\times\EG:p\mapsto(u(p),\phi(p))$
is equivariant and descends to a function
$f:X\to M\times_\G\EG$ that satisfies
$$
    f_*[X] = [u].
$$
To see this, consider the maps
$
    \phi^\G:X\to P\times_\G\EG
$
and
$
    \pi^\G:P\times_\G\EG\to X, 
$
given by 
$
    \phi^\G(\pi(p))=[p,\phi(p)]
$
and
$
    \pi^\G([p,e])=\pi(p).
$
Since $\EG$ is contractible, $\phi^\G$ 
is a homotopy inverse of $\pi^\G$. 
Hence $\phi^\G_*$ is the inverse of $\pi^\G_*$. 
Moreover, 
$
    f=u^\G\circ\phi^\G,
$
where $u^\G:P\times_\G\EG\to M\times_\G\EG$
is given by $u^\G([p,e])=[u(p),e]$.
Hence $f_*[X]=u^\G_*\phi^\G_*[X]=[u]$. 
\end{remark}

\begin{proposition}\label{prop:2}
Let $M$ be a finite dimensional smooth manifold
and $\G$ be a compact Lie group 
which acts smoothly on $M$.

\smallskip
\noindent{\bf (i)} 
For every $2$-dimensional equivariant homology class
$B\in H^\G_2(M;\Z)$ there exists a compact oriented
Riemann surface $\Sigma$, a principal bundle $P\to\Sigma$,
and an equivariant map $u:P\to M$, such that $[u]=B$.

\smallskip
\noindent{\bf (ii)}
Suppose that $\G$ is connected.
Let $P\to\Sigma$ and $P'\to\Sigma$ be principal 
$\G$-bundles over $\Sigma$, 
and $u:P\to M$, $u':P'\to M$ be equivariant
maps such that $[u]=[u']\in H^\G_2(M;\Z)$.
Then $P$ is isomorphic to $P'$. 
\end{proposition}

\begin{proof}
Given $B$, choose a compact oriented Riemann surface 
$\Sigma$ and a map 
$
     f:\Sigma\to M\times_\G\EG
$
such that $f_*[\Sigma]=B$. Note that $M\times\EG$
is a principal $\G$-bundle over $M\times_\G\EG$ and 
denote by $P\to\Sigma$ the pullback bundle of $f$.
Thus 
$$
     P = \left\{(z,x,e)\in\Sigma\times M\times\EG\,|\,
        [x,e]=f(z)\right\}.
$$
There are two equivariant maps 
$u:P\to M$ and $\phi:P\to\EG$, given by,
$$
     u(z,x,e) = x,\qquad
     \phi(z,x,e) = e.
$$
By definition, the map 
$(u,\phi):P\to M\times\EG$
descends to $f$.  Hence,
by Remark~\ref{rmk:P},
$[u]=f_*[\Sigma]=B$. 
This proves~(i).

To prove~(ii), choose two equivariant maps 
$\phi:P\to\EG$ and $\phi':P'\to\EG$.
Define $f,f':\Sigma\to M\times_\G\EG$ as the maps
induced by $(u,\phi):P\to M\times\EG$
and $(u',\phi'):P'\to M\times\EG$. 
Then, by Remark~\ref{rmk:P},
$$
     f_*[\Sigma] = [u] = [u'] = {f'}_*[\Sigma].
$$
Consider the induced maps 
$
     \bar\phi:\Sigma\to\BG
$ 
and 
$
     \bar\phi':\Sigma\to\BG.
$  
They can be expressed in the
form
$
    \bar\phi=\pi\circ f
$
and
$
    \bar\phi'=\pi\circ f',
$
where $\pi:M\times_\G\EG\to\BG$ denotes the obvious
projection. Hence $\bar\phi$ and $\bar\phi'$ are homologous, 
i.e. $\bar\phi_*[\Sigma]=\bar\phi'_*[\Sigma]$. 
Since $\G$ is connected, $\BG$ is simply connected.  
Hence two maps $\Sigma\to\BG$ are homologous 
if and only they are homotopic. (To see this note that
every map from $\Sigma$ to a simply connected 
space factors, up to homotopy, through a map 
of degree $1$ from $\Sigma$ to $S^2$.)
This shows that our maps $\bar\phi$ and $\bar\phi'$ 
are homotopic. 
Hence $P$ and $P'$ are isomorphic. 
This proves the proposition. 
\end{proof}

Assertion~(ii) in Proposition~\ref{prop:2} can be restated 
as follows.  An equivariant homology class $B\in H^\G_2(M;\Z)$
descends to a homology class $b\in H_2(\BG;\Z)$ and the 
latter determines an isomorphism class of principal
$\G$-bundles $P\to\Sigma$ 
(over any orientable Riemann surface). 

\smallbreak

The deRham model of equivariant cohomology
is defined as follows.   Let $\Om^*_\G(M)$
denote the space of equivariant 
polynomials from $\g$ to $\Om^*(M)$.
To be more explicit, choose a basis 
$e_1,\dots,e_m$ of $\g$ and write
$
    \xi=\sum_{i=1}^m\xi^ie_i\in\g.
$
Then any $\alpha\in\Om^k_\G(M)$ can be written 
in the form
$$
    \alpha(\xi) = \sum_I\xi^I\alpha_I,
$$
where $I=(i_1,\dots,i_m)$, $\xi^I=(\xi^1)^{i_1}\cdots(\xi^m)^{i_m}$,
and $\alpha_I\in\Om^{k-2|I|}(M)$.
The equivariance of the function
$\alpha:\g\to\Om^*(M)$ can be expressed in the form
$$
    D\alpha(\xi)[\xi,\eta] = \Ll_{X_\eta}\alpha(\xi)
$$
for $\xi,\eta\in\g$. Here the linear operator 
$D\alpha(\xi):\g\to\Om^*(M)$ denotes
the differential of the function $\g\to\Om^*(M):\xi\mapsto\alpha(\xi)$
at the point $\xi$. 
The differential 
$
    d_\G:\Om^k_\G(M)\to\Om^{k+1}_\G(M)
$
is given by
$$
    (d_\G\alpha)(\xi) 
    := d(\alpha(\xi)) + \i(X_\xi)\alpha(\xi)
    = \sum_I\xi^I(d\alpha_I+\i(X_\xi)\alpha_I).
$$
Cartan's formula asserts that $d_\G\circ d_\G=0$. 
The equivariant version of de\-Rham's theorem
asserts that, for every smooth manifold $M$
with a smooth $\G$-action, there is a natural
isomorphism~\cite{KIRWAN}
$$
    H^k_\G(M;\R) \cong
    \frac{\ker d_\G:\Om^k_\G(M)\to\Om^{k+1}_\G(M)}
    {\im d_\G:\Om^{k-1}_\G(M)\to\Om^k_\G(M)}.
$$
Next we describe the pairing between an
equivariant cohomology class, represented by 
a $\G$-closed $k$-form $\alpha\in\Om_\G^k(M)$
and an equivariant homology class, represented 
by an equivariant map $u:P\to M$ defined
on the total space of a principal $\G$-bundle 
$\pi:P\to X$ over a compact oriented smooth 
$k$-manifold $X$.  An explicit formula for this 
pairing relies on a $\G$-connection $A\in\Aa(P)$
and the covariant derivative of $u$ determined by $A$.
This covariant derivative is defined as follows.
Think of $u$ as a section of the 
associated bundle 
$$
     \Tilde{M}=P\times_\G M\to X
$$
with fibres diffeomorphic to $M$.
The connection $A$ on $P$ determines a 
symplectic connection on this bundle.  
More precisely, the tangent space 
of $\Tilde{M}$ at $[p,x]$ is the quotient
$$
     T_{[p,x]}\Tilde{M}
     = T_pP\times T_xM/\{(p\xi,-X_\xi(x))\,|\,\xi\in\g\},
$$
the vertical space consists of equivalence classes 
of the form $[0,w]$ with $w\in T_xM$,
and the horizontal space 
consists of those equivalence classes 
$[v,w]$, where $v\in T_pP$ and $w\in T_xM$
satisfy $w+X_{A_p(v)}(x)=0$. 
In the terminology of~\cite[Chapter~6]{MS2}
the connection form on $\Tilde{M}$ 
is induced by the $2$-form
$\om-d\inner{\mu}{A}$ on $P\times M$,
whenever $(M,\om,\mu)$ is a symplectic manifold
with a Hamiltonian group action.
The covariant derivative 
of a section $u:P\to M$ with respect to the 
connection $A$ is the $1$-form 
$
     d_Au : TP\to u^*TM
$
given by 
\begin{equation}\label{eq:dAu}
    d_Au(p)v = du(p)v + X_{A_p(v)}(u(p))
\end{equation}
(the vertical part of the vector 
$[v,du(p)v]\in T_{[p,u(p)]}P\times_\G M$). 
This $1$-form is obviously $\G$-equivariant,
and it satisfies
$
     d_Au(p)p\xi = 0
$
for every $\xi\in\g$.  Hence $d_Au$ descends to a $1$-form
on $X$ with values in the bundle $u^*TM/\G$.  

Given a basis $e_1,\dots,e_m$ of $\g$ 
and an equivariant $k$-form 
$\alpha=\sum_I\xi^I\alpha_I\in\Om^k_\G(M)$ 
as above, we define $\alpha(u,A)\in\Om^k(P)$ by
$$
     \alpha(u,A) 
     = ((d_Au)^*\alpha)(F_A)
     = \sum_I\om^I\wedge(d_Au)^*\alpha_I.
$$
Here $F_A=\sum_i\om^ie_i$ and 
$\om^I=(\om^1)^{i_1}\wedge\cdots\wedge(\om^m)^{i_m}$.
Since $\alpha(u,A)\in\Om^k(P)$ is equivariant 
and horizontal (see Proposition~\ref{prop:eco} below)
it descends to a $k$-form on $X$, 
still denoted by $\alpha(u,A)$.  
The pairing between the equivariant cohomology class 
of $\alpha$ and the equivariant homology class of $u$ 
is given by 
$$
     \inner{[\alpha]}{[u]}
     := \int_X\alpha(u,A).
$$
That this number is well defined and depends only 
on the equivariant cohomology class of $\alpha$
and on the homotopy class of the pair $(u,A)$
is the content of the following proposition.
 
\begin{proposition}\label{prop:eco}
Let $M$ be a smooth $\G$-manifold and
$\pi:P\to X$ be a principal $\G$-bundle over 
a compact smooth manifold $X$.  
Let $\alpha\in\Om^\ell_\G(M)$
and $\beta\in\Om^{\ell+1}_\G(M)$.

\smallskip
\noindent{\bf (i)}
Let $A\in\Aa(P)$ and suppose that 
$u:P\to M$ is an equivariant smooth function.
Then $\alpha(u,A)\in\Om^\ell(P)$ is equivariant and
horizontal.  

\smallskip
\noindent{\bf (ii)}
If $d_\G\alpha=\beta$ then $d\alpha(u,A)=\beta(u,A)$.

\smallskip
\noindent{\bf (iii)}
Let $A_0,A_1\in\Aa(P)$ and suppose that 
$u_0,u_1:P\to M$ are equivariantly homotopic.  Then 
$\alpha(u_1,A_1)-\alpha(u_0,A_0)$ is an exact
$\ell$-form on $X$. 
\end{proposition} 
 
\begin{proof}
The form $\alpha(u,A)$ is obviously horizontal.
We prove that it is equivariant. 
Denote by $c^k_{ij}$ the structure constants of $\g$.
This means that
$$
     [e_i,e_j]=\sum_k c^k_{ij}e_k.
$$
Then the equivariance of $\alpha$ can be expressed in the form
$$
     \sum_{i,k} c^k_{ij}\xi^i\alpha_k(\xi) 
     = \Ll_{X_{e_j}}\alpha(\xi),
$$
where $\alpha_k:=\p_k\alpha:\g\to\Om^*(M)$. 
Hence, with 
$
     F_A 
     = \om 
     = \sum_i\om^ie_i\in\Om^2(P,\g),
$ 
it follows that
\begin{equation}\label{eq:equi-alpha}
     \sum_{i,k} c^k_{ij}\om^i\wedge\alpha_k(u,A) 
     = (\Ll_{X_{e_j}}\alpha)(u,A).
\end{equation}
Moreover, with $A=\sum_ia^ie_i\in\Om^1(P,\g)$,
the Bianchi identity $d_AF_A=0$
takes the form 
\begin{equation}\label{eq:bianchi}
      d\om^k = \sum_{i,j} c^k_{ij} \om^i\wedge a^j.
\end{equation}
For $j=1,\dots,m$ denote by $v_j\in\Vect(P)$ 
the vector field $p\mapsto pe_j$. 
Then $a^i(v_j)=\delta^i_j$.
Hence, by~(\ref{eq:bianchi}),
$$
     \Ll_{v_j}\om^k = \i(v_j)d\om^k = \sum_i c^k_{ij}\om^i,
$$
and hence, by~(\ref{eq:equi-alpha}),
$$
      \sum_{k}(\Ll_{v_j}\om^k)\wedge\alpha_k(u,A)  
      = (\Ll_{X_{e_j}}\alpha)(u,A) 
      = - \sum_I\om^I\wedge \Ll_{v_j}((d_Au)^*\alpha_I) .
$$
Here we have used the identity
$
     (d_Au)^*\Ll_{X_{e_j}}\alpha_I
     = - \Ll_{v_j}((d_Au)^*\alpha_I).
$
It follows that 
$$
     \Ll_{v_j}\alpha(u,A) 
     = \sum_{k}(\Ll_{v_j}\om^k)\wedge\alpha_k(u,A) 
       + \sum_I\om^I\wedge \Ll_{v_j}((d_Au)^*\alpha_I) 
     = 0
$$
and hence $\alpha(u,A)$ is equivariant.
This proves~(i).

The proof of~(ii) relies on the following identity,
for $\alpha\in\Om^k(M)$,
\begin{eqnarray}\label{eq:ddAu}
&&
    d((d_Au)^*\alpha) - (d_Au)^*d\alpha
    \nonumber  \\
&&\quad
    = \sum_i
      \om^i\wedge (d_Au)^*\i(X_{e_i})\alpha
      - \sum_ia^i\wedge(d_Au)^*\Ll_{X_{e_i}}\alpha.
\end{eqnarray}
For $k=0,1$ the proof of~(\ref{eq:ddAu})
is a computation using
$
    \om^k=da^k+\sum_{i<j}c^k_{ij}a^i\wedge a^j
$
and
$$
    \i(X_{e_j})\Ll_{X_{e_i}}\alpha 
    - \i(X_{e_i})\Ll_{X_{e_j}}\alpha
    - d\alpha(X_{e_i},X_{e_j})
    = \alpha([X_{e_i},X_{e_j}]).
$$
For general $k$, (\ref{eq:ddAu}) 
follows easily by induction.  With this 
understood, we obtain,
\begin{eqnarray*}
     d\alpha(u,A)
&= &
     \sum_kd\om^k\wedge\alpha_k(u,A)
     + \sum_{I}\om^I\wedge d((d_Au)^*\alpha_I) \\
&= &
     \sum_{I,j}\om^I\wedge a^j
     \wedge(d_Au)^*\Ll_{X_{e_j}}\alpha_I 
     + \sum_{I}\om^I\wedge d((d_Au)^*\alpha_I) \\
&= &
     \sum_{I,j}\om^I\wedge\om^j\wedge (d_Au)^*\i(X_{e_j})\alpha_I
     + \sum_{I}\om^I\wedge(d_Au)^*d\alpha_I  \\
&= &
     \beta(u,A)
\end{eqnarray*}
Here the second identity follows from~(\ref{eq:equi-alpha})
and~(\ref{eq:bianchi}), the third identity from~(\ref{eq:ddAu}),
and the last identity from $d_\G\alpha=\beta$, i.e.
$
     \sum_I \xi^I(d\alpha_I+\i(X_\xi)\alpha_I)
     = \sum_I \xi^I\beta_I.
$ 
This proves~(ii).  

We prove~(iii).  Let $\R\to\Aa(P):t\mapsto A_t$
be a smooth family of connections and 
$\R\times P\to M:(t,p)\mapsto u_t(p)$ 
be a smooth family of equivariant functions.
Think of the path $t\mapsto A_t$ as a connection $\Tilde{A}$
on the bundle $\Tilde{P}=\R\times P$ over $\Tilde{X}=\R\times X$
and of the path $t\mapsto u_t$ as a function 
$\Tilde{u}:\Tilde{P}\to M$.  Given a $\G$-closed $\ell$-form
$\alpha\in\Om^\ell_\G(M)$ write
$$
    \alpha(\Tilde{u},\Tilde{A})
    = \Tilde{\alpha} = \alpha_t + \beta_t\wedge dt
    \in \Om^{\ell}(\Tilde{P}),
$$
where $\alpha_t=\alpha(u_t,A_t)\in\Om^\ell(P)$
and $\beta_t\in\Om^{\ell-1}(P)$. 
By~(i), $\Tilde{\alpha}$ is horizontal, equivariant, and closed.
Hence $\alpha_t$ and $\beta_t$ are equivariant and horizontal,
$\alpha_t$ is closed, and  
$\p_t\alpha_t=d\beta_t$ 
for every $t$. Hence 
$$
   \alpha(u_1,A_1)  - \alpha(u_0,A_0)
   = d\int_0^1\beta_t\,dt.
$$
Since $\beta_t$ descends to $X$ for every $t$,
this proves the proposition. 
\end{proof}


\subsection{J-holomorphic curves}\label{sec:jhol}

Suppose that $(M,\om,\mu)$ is a symplectic 
manifold with a Hamiltonian group action.
Denote by $\Jj(M,\om,\mu)$
the space of all almost complex structures $J$ on $TM$
that are invariant under the $\G$-action 
and are compatible with $\om$, i.e. $\om(\cdot,J\cdot)$
is a Riemannian metric on $M$. 
It follows from Proposition~2.50 in~\cite{MS2}
that the space $\Jj(M,\om,\mu)$ is nonempty 
and contractible.  Namely, there is a natural
homotopy equivalence from the (contractible)
space of $\G$-invariant Riemannian metrics on $M$
to the space $\Jj(M,\om,\mu)$.

An almost complex structure $J\in\Jj(M,\om,\mu)$ 
determines an almost complex structure 
on the quotient 
$$
     M\dslash\G = \mu^{-1}(0)/\G.
$$  
The tangent space of this quotient is given by 
$$
     T_{[x]}M\dslash\G = \ker d\mu(x)\cap\ker{L_x}^*
$$
and the identity
$$
     d\mu(x)J={L_x}^*
$$
shows that this space is invariant under $J$. 
Hence a map $u:\C\to\mu^{-1}(0)$ represents
a $J$-holomorphic curve in $M\dslash\G$ if and only if there
exist functions $\Phi,\Psi:\C\to\g$ such that
\begin{equation}\label{eq:jhol-loc}
    \p_su+X_\Phi(u) + J(\p_tu+X_\Psi(u)) = 0.
\end{equation}
Here we denote by $s+it$ the coordinate on $\C$. 
This equation implies that the vectors
$\p_su+X_\Phi(u)$ and $\p_tu+X_\Psi(u)$
are the unique harmonic representatives
of the derivatives with respect to $s$ and $t$.
Thus they are uniquely determined
by the equations
$$
    {L_u}^*L_u\Phi + {L_u}^*\p_su = 0,\qquad
    {L_u}^*L_u\Psi + {L_u}^*\p_tu = 0.      
$$
There is a natural gauge group of maps $g:\C\to\G$.
It acts on triples $(u,\Phi,\Psi)$ by 
$$
    g^*(u,\Phi,\Psi)
    = (g^{-1}u,
       g^{-1}\p_sg+g^{-1}\Phi g,
       g^{-1}\p_tg+g^{-1}\Psi g).
$$
This action preserves the space of solutions 
of~(\ref{eq:jhol-loc}).
Note that the action on the space of $1$-forms 
$\Phi\,ds+\Psi\,dt$ coincides with the action 
of the gauge group on the space of connections. 

The global version of equation~(\ref{eq:jhol-loc}) 
involves principal bundles.  
Let $\Sigma$ be a compact oriented Riemann surface 
with a fixed complex structure $J_\Sigma$. 
Then a smooth function
$\Sigma\to M\dslash\G$ need not lift 
to a smooth function into the ambient space $M$.  
However, it does lift to 
an equivariant function from the total
space of a principal $\G$-bundle 
$\pi:P\to\Sigma$ into $M$.  
Hence let $(u,A)$ be a pair consisting
of an equivariant smooth function 
$u:P\to M$ and a connection $A\in\Aa(P)$.
Recall that $d_Au:TP\to u^*TM$ is defined by 
$$
     d_Au = du + L_uA.
$$
Think of $d_Au$ as a $1$-form on $\Sigma$ 
with values in the bundle $u^*TM/\G$.
This is a complex vector bundle
and we shall denote by 
$\bar\p_{J,A}(u)\in\Om^{0,1}(\Sigma,u^*TM/\G)$
the complex anti-linear part 
of the $1$-form $d_Au$.  Thus
$$
    \bar\p_{J,A}(u) 
        = \frac12(d_Au+J\circ d_Au\circ J_\Sigma)
        \in\Om^{0,1}(\Sigma,u^*TM/\G).
$$
In local coordinates this is the left hand side
of~(\ref{eq:jhol-loc}).  To make sense of this expression
in the global form note that $J_\Sigma$ acts on the tangent 
space of $\Sigma$ but not on that of $P$.  
The linear map $d_Au(p)\circ J_\Sigma:T_pP\to T_{u(p)}M$
is defined as follows.  Project a vector $v\in T_pP$
onto $T_{\pi(p)}\Sigma$ and then apply $J_\Sigma$. 
Now lift $J_\Sigma d\pi(p)v\in T_{\pi(p)}\Sigma$ to
a vector in $T_pP$ and apply $d_Au(p)$.  Since 
$d_Au$ is horizontal, the resulting vector in  $T_{u(p)}M$
is independent of the choice of the lift. 
To understand the $(0,1)$-form $\bar\p_{J,A}(u)$ 
in a different way consider the fibre bundle 
$\Tilde{M}=P\times_\G M\to\Sigma$,
with fibres diffeomorphic to $M$,
and define $\Tilde{u}:\Sigma\to\Tilde{M}$
by $\Tilde{u}(\pi(p))=[p,u(p)]$. 
Now $J_\Sigma$, $J$, and $A$ determine an almost complex
structure $\Tilde{J}_A$ on $\Tilde{M}$,\label{tildeJA} 
and $\Tilde{u}$ is a $\Tilde{J}_A$-holomorphic curve 
if and only if $\bar\p_{J,A}(u)=0$.
In any case, the global form of~(\ref{eq:jhol-loc}) is 
\begin{equation}\label{eq:jhol}
     \bar\p_{J,A}(u)=0,\qquad
     \mu(u)=0.
\end{equation}
Here $A\in\Aa(P)$ is the pullback under $u$ of the 
connection on the principle bundle $\mu^{-1}(0)\to M\dslash\G$ 
that is determined by the metric $\om(\cdot,J\cdot)$. 
Fix an equivariant homology class 
$B\in H_2^\G(M;\Z)$ and denote by
$$
     \Mm_{B,\Sigma}^0(M,\mu;J)
     = \left\{(u,A)\,|\,(\ref{eq:jhol}),\,[u]=B
       \right\}/\Gg(P)
$$
the moduli space of gauge equivalence classes of solutions
of~(\ref{eq:jhol}) that represent the class $B$. 
If $B$ is the image of a class 
$\bar B\in H_2(M\dslash\G;\Z)$
under the natural homomorphisms
$
     H_2(M\dslash\G;\Z)
     \to H_2^\G(M;\Z)
$
then the solutions of~(\ref{eq:jhol}) 
correspond to $J$-holomorphic curves
in the quotient $M\dslash\G$ representing the 
class $\bar B$. 

\begin{remark}\label{rmk:tau}\rm
{\bf (i)} If $\tau\in\g$ is a central element
then we can replace the moment map $\mu$ 
by $\mu-\tau$.  The solutions of 
\begin{equation}\label{eq:jhol-tau}
     \bar\p_{J,A}(u)=0,\qquad
     \mu(u)=\tau
\end{equation}
correspond to $J$-holomorphic curves in the 
quotient $M\dslash\G(\tau)$.

\smallskip
\noindent{\bf (ii)}
Let $\Oo\subset\g$ be an orbit under the adjoint action
and consider the product $M\times\Oo$ 
with the moment map $\mu_\Oo:M\times\Oo\to\g$
given by $\mu_\Oo(x,\tau)=\mu(x)-\tau$. 
Then~(\ref{eq:jhol}) takes the form
\begin{equation}\label{eq:jhol-Oo}
     \bar\p_{J,A}(u)=0,\qquad
     \bar\p_A(\tau)=0,\qquad
     \mu(u)=\tau,
\end{equation}
where $u:P\to M$ and $\tau:P\to\Oo$ are equivariant maps and
$$
     \bar\p_A(\tau) 
     = \frac12(d_A\tau-J_\tau\circ d_A\tau\circ J_\Sigma),\qquad
     d_A\tau = d\tau + [A,\tau].
$$
The solutions of~(\ref{eq:jhol-Oo}) again correspond to
$J$-holomorphic curves in the quotient $M\dslash\G(\tau)$.
Note that~(\ref{eq:jhol-Oo}) is equivalent to~(\ref{eq:jhol-tau})
whenever $\Oo$ is a single point (necessarily
contained in the centre of $\g$). In local holomorphic
coordinates~(\ref{eq:jhol-Oo}) has the form
\begin{equation}\label{eq:jhol-loc-Oo}
\begin{array}{rcl}
    \p_su+X_\Phi(u) + J(\p_tu+X_\Psi(u)) &= &0,  \\
    \p_s\tau+[\Phi,\tau] - J_\tau(\p_t\tau+[\Psi,\tau]) &= & 0, \\
    \mu(u)-\tau &= & 0.
\end{array}
\end{equation}
As before, these equations imply that the pairs 
$(\p_su+X_\Phi(u),\p_s\tau+[\Phi,\tau])$ and
$(\p_tu+X_\Psi(u),\p_t\tau+[\Psi,\tau])$
are the unique harmonic representatives
of the derivatives with respect to $s$ and $t$.
Thus, by Remark~\ref{rmk:quotient}~(ii), 
the function $\Phi:\C\to\g$ 
is determined by the equation
$$
    {L_u}^*L_u\Phi+A_\tau\Phi + {L_u}^*\p_su + J_\tau\p_s\tau = 0,
$$
and similarly for $\Psi$.  Note that in this local form
the quadruple $(u,\tau,\Phi,\Psi)$ can be gauge 
transformed to one where $\tau$ is constant. 
The second equation in~(\ref{eq:jhol-loc-Oo})
then takes the form 
$
    [\Phi,\tau]=J_\tau[\Psi,\tau]
$
and can be viewed as a constraint on the connection $A$.
The gauge group now consists of maps $\C\to\G_\tau$. 

\smallskip
\noindent{\bf (iii)}
If one removes the assumption that $\G$ acts freely on 
$\mu^{-1}(0)$, then $M\dslash\G$ is an orbifold
and~(\ref{eq:jhol-Oo}) describes $J$-holomorphic curves 
in this space.
Since every symplectic orbifold can be expressed in this 
form (cf.~\cite{MARTIN3}), one might be tempted 
to use~(\ref{eq:jhol}) to give a rigorous definition 
of the Gromov--Witten invariants of orbifolds. 
\end{remark}


\section{Invariants of Hamiltonian group actions}\label{sec:new}


\subsection{An action functional}

Let $(M,\om,\mu)$ be a symplectic manifold 
with a Hamiltonian $\G$-action and
$\pi:P\to\Sigma$ be a principle $\G$-bundle
over a compact Riemann surface $(\Sigma,J_\Sigma)$.
Denote by $\Cinf_\G(P,M)$ the space of equivariant 
smooth functions $u:P\to M$ and consider the action
functional $E:\Cinf_\G(P,M)\times\Aa(P)\to\R$, 
defined by 
$$
     E(u,A) 
     = \frac12\int_\Sigma
       \left(
         \left|d_Au\right|^2
         + \left|F_A\right|^2
         + \left|\mu(u)\right|^2
       \right)
       \dvol_\Sigma.
$$
This functional is invariant under the action
of the gauge group $\Gg(P)$.
The Euler equations have the form
\begin{equation}\label{eq:euler}
     {\Nabla{A}}^*d_Au+d\mu(u)^*\mu(u)=0,\qquad
     {d_A}^*F_A + {L_u}^*d_Au = 0.
\end{equation}
Here 
$
     \Nabla{A}:\Cinf(\Sigma,u^*TM/\G)
     \to\Om^1(\Sigma,u^*TM/\G)
$
denotes the covariant derivative operator 
induced by $A$ and by the Levi-Civita connection 
$\nabla$ of the metric $\om(\cdot,J\cdot)$ on $M$.
It is defined by 
\begin{equation}\label{eq:nablaA}
    \Nabla{A}\xi = \nabla\xi + \Nabla{\xi}X_A
\end{equation}
for $\xi\in\Cinf(\Sigma,u^*TM/\G)$.
Think of $\xi$ as an equivariant function from
$P$ to $u^*TM$.  Then $\Nabla{A}\xi$ is a 
$1$-form on $P$ with values in $u^*TM$.
This form is obviously equivariant and,
since $\Nabla{p\eta}\xi(p)+\Nabla{\xi(p)}X_\eta(u(p))=0$,
it is horizontal.  Hence it descends to a 
$1$-form on $\Sigma$ with values in $u^*TM/\G$,
still denoted by $\Nabla{A}\xi$.
The symbol ${\Nabla{A}}^*$ in~(\ref{eq:euler})
denotes the $L^2$-adjoint of $\Nabla{A}$.

There are first order equations that give rise to 
special solutions of~(\ref{eq:euler}). They have 
the form
\begin{equation}\label{eq:new}
     \bar\p_{J,A}(u) = 0,\qquad
     *F_A+\mu(u) = 0.
\end{equation}
Here $*$ denotes the Hodge $*$-operator on $\Sigma$
and so these equations depend explicitly 
on the metric on $\Sigma$.
The study of their solutions
is the main purpose of this paper.  
The next proposition  shows that the solutions
of~(\ref{eq:new}), if they exist, are the absolute
minima of $E$ and hence also solve the Euler 
equations.  The moment map condition asserts
that the polynomial
$$
      \g\to\Om^*(M):\xi\mapsto \om-\inner{\mu}{\xi}
$$
is $\G$-closed and hence defines an equivariant
cohomology class which we denote by 
$[\om-\mu]\in H^2_\G(M;\R)$. 
Hence, by Proposition~\ref{prop:eco},
the last term in~(\ref{eq:energy}) below
is a topological invariant. 

\begin{proposition}\label{prop:energy}
For every $A\in\Aa(P)$ and every $u\in\Cinf_\G(P,M)$,
\begin{equation}\label{eq:energy}
     E(u,A) 
     = \int_\Sigma
       \left(\left|\bar\p_{J,A}(u)\right|^2
       + \frac12\left|*F_A+\mu(u)\right|^2\right)\dvol_\Sigma 
       + \inner{[\om-\mu]}{[u]},
\end{equation}
where
$$
     \inner{[\om-\mu]}{[u]}
     = \int_\Sigma \left((d_Au)^*\om-\inner{\mu(u)}{F_A}\right)
     = \int_\Sigma \left(u^*\om-d\inner{\mu(u)}{A}\right).
$$
\end{proposition}

\begin{proof}  Choose a holomorphic coordinate chart
$\phi:U\to\Sigma$, where $U\subset\C$ is an open set,
and let $\tilde\phi:U\to P$ be a lift of $\phi$,
that is $\pi\circ\tilde\phi=\phi$.   
Then the function $u$ and the connection $A$
are in local coordinates given 
by $u^\loc=u\circ\tilde\phi$
and $A^\loc={\tilde\phi\,}^*A=\Phi\,ds+\Psi\,dt$
where $\Phi,\Psi:U\to\g$.  The pullback volume form 
on $U$ is
$
     \lambda^2\,ds\wedge dt
$
for some function $\lambda:U\to(0,\infty)$
and the metric is $\lambda^2(ds^2+dt^2)$.
Hence
$$
     {\tilde\phi\,}^*F_A
         = \left(\p_s\Psi-\p_t\Phi+[\Phi,\Psi]\right) 
           \,ds\wedge dt,
$$
$$
     {\tilde\phi\,}^*d_Au
         = \left(\p_su^\loc+X_\Phi\right)\,ds 
     + \left(\p_tu^\loc+X_\Psi\right)\,dt,
$$
$$
     {\tilde\phi\,}^*\bar\p_{J,A}(u)
     = \frac12(\xi\,ds -J\xi\,dt),\quad
     \xi = \p_su^\loc+X_\Phi+J(\p_tu^\loc+X_\Psi). 
$$
Here $X_\Phi$, $X_\Psi$, and $J$ are evaluated at $u^\loc$. 
In the following we shall drop the superscript ``$\loc$''.
Then the equations~(\ref{eq:new}) have the form
\begin{equation}\label{eq:new-loc}
\begin{array}{rcl}
    \p_su+X_\Phi(u) + J(\p_tu+X_\Psi(u)) &= &0,  \\
    \p_s\Psi-\p_t\Phi+[\Phi,\Psi] + \lambda^2\mu(u) &= & 0,
\end{array}
\end{equation}
The pullback of the energy integrand under $\phi:U\to\Sigma$
is given by 
\begin{eqnarray*}
    e: 
&= &
    \frac12\left|\p_su+X_\Phi\right|^2
        + \frac12\left|\p_tu+X_\Psi\right|^2  \\
&&
    + \,\frac{\lambda^{-2}}{2}\left|\p_s\Psi-\p_t\Phi+[\Phi,\Psi]\right|^2
    + \frac{\lambda^2}{2}\left|\mu(u)\right|^2  \\
&= &
    \frac12\left|\p_su+X_\Phi+J(\p_tu+X_\Psi)\right|^2 \\
&&
    +\,\frac{\lambda^2}{2}
        \left|\lambda^{-2}\left(\p_s\Psi-\p_t\Phi+[\Phi,\Psi]\right) 
         + \mu(u)\right|^2   \\
&&
    +\,\om(\p_su+X_\Phi,\p_tu+X_\Psi)  
    - \inner{\p_s\Psi-\p_t\Phi+[\Phi,\Psi]}{\mu(u)}.
\end{eqnarray*}
This proves~(\ref{eq:energy}).
The identity
$
     (d_Au)^*\om-\inner{\mu(u)}{F_A}
     = u^*\om-d\inner{\mu(u)}{A}
$
can in local coordinates be expressed 
in the form
\begin{eqnarray*}    
&&
    \om(\p_su+X_\Phi,\p_tu+X_\Psi)  
    - \inner{\p_s\Psi-\p_t\Phi+[\Phi,\Psi]}{\mu(u)}  \\
&&\qquad=
    \om(\p_su,\p_tu)
    - \p_s\inner{\mu(u)}{\Psi}
    + \p_t\inner{\mu(u)}{\Phi}.
\end{eqnarray*}
This follows directly from the definitions and the fact that
$\om(X_\Phi,X_\Psi)=\inner{\mu}{[\Phi,\Psi]}$.
This proves the proposition.
\end{proof}

%
%


\subsection{Symplectic reduction}\label{subsec:red}

Denote $\Cinf_\G(P,M;B)=\{u\in\Cinf_\G(P,M)\,|\,[u]=B\}$
and consider the space\label{Bb}
$$
     \Bb = \Cinf_\G(P,M;B) \times \Aa(P).
$$
This space carries a natural symplectic form.
To see this note that the tangent space
of $\Bb$ at $(u,A)$ is
$$
     T_{(u,A)}\Bb 
         = \Cinf(\Sigma,u^*TM/\G)\oplus\Om^1(\Sigma,\g_P),
$$
where $\Cinf(\Sigma,u^*TM/\G)=\Cinf_\G(P,u^*TM)$
can be thought of as the space of $\G$-equivariant
sections of the bundle $u^*TM\to P$ and 
$\Om^1(\Sigma,\g_P)=\Om^1_\ad(P,\g)$ as the space 
of equivariant and horizontal Lie algebra 
valued $1$-forms on $P$.  The symplectic form on $T_{(u,A)}\Bb$
is given by 
\begin{equation}\label{eq:Om}
     \Om((\xi,\alpha),(\xi',\alpha'))
     = \int_\Sigma\om(\xi,\xi')\dvol_\Sigma
           + \int_\Sigma\winner{\alpha}{\alpha'}
\end{equation}
for $\xi,\xi'\in\Cinf(\Sigma,u^*TM/\G)$ and 
$\alpha,\alpha'\in\Om^1(\Sigma,\g_P)$
(see Section~\ref{sec:gauge}).
Consider the group 
$\Tilde{\Gg}=\Tilde{\Gg}(P)$
of all automorphisms of $P$ that descend to
Hamiltonian symplectomorphisms of $\Sigma$.
There is an exact sequence
$$
    1
    \TO\Gg
    \TO\Tilde{\Gg}
    \TO\Ham(\Sigma,\dvol_\Sigma)
    \TO 1
$$
and the Lie algebra of $\Tilde{\Gg}$ consists of all
equivariant vector fields $v\in\Vect_\G(P)$ 
such that the $1$-form $\i(\pi_*v)\dvol_\Sigma\in\Om^1(\Sigma)$
is exact.  Thus, for every $v\in\Lie(\Tilde{\Gg})$,
there exists a unique function $h_v:\Sigma\to\R$ such that
\begin{equation}\label{eq:h}
      \i(\pi_*v)\dvol_\Sigma=dh_v,\qquad
      \int_\Sigma h_v\dvol_\Sigma=0.
\end{equation}
The group $\Tilde{\Gg}$ acts on $\Bb$ by 
$
     (u,A)\mapsto(u\circ f^{-1},f_*A)
$
for $f\in\Tilde{\Gg}$ and the infinitesimal action is given by 
the vector fields
$$
     \Bb\to T\Bb:
     (u,A) \mapsto (-du\circ v,-\Ll_vA)
$$
for $v\in\Lie(\Tilde{\Gg})$.
It follows from the work of Atiyah--Bott~\cite{AB}
and Donaldson~\cite{DONALDSON6} that this action is Hamiltonian. 
Define $\Tilde{\mu}:\Bb\to\Lie(\Tilde{\Gg})^*$ by 
\begin{equation}\label{eq:Moment}
     \inner{\Tilde{\mu}(u,A)}{v}
     = \int_\Sigma \biggl(
        \inner{*F_A+\mu(u)}{A(v)}\dvol_\Sigma
        - h_v\left(u^*\om-d\inner{\mu(u)}{A}\right)
        \biggr).
\end{equation}
Here $\inner{*F_A+\mu(u)}{A(v)}\in\Om^0(P)$ and
and $u^*\om-d\inner{\mu(u)}{A}\in\Om^2(P)$.  
But they both descend to $\Sigma$ and we do not
distinguish the descendents in notation from the 
original forms on $P$. 

\begin{proposition}\label{prop:Moment}
The function $\Tilde{\mu}:\Bb\to\Lie(\Tilde{\Gg})^*$
is a moment map for the action 
of $\Tilde{\Gg}$ on $\Bb$. 
\end{proposition}

\begin{proof}
Let $\R\to\Bb:t\mapsto(u_t,A_t)$ be any smooth
path in $\Bb$ and denote its derivative by 
$
     (\xi_t,\alpha_t)\in T_{(u_t,A_t)}\Bb.
$
Then, by Cartan's formula,
$$
     \frac{d}{dt}({u_t}^*\om-d\inner{\mu(u_t)}{A_t}) 
     = d\sigma_t,\quad
     \sigma_t 
     = \om(\xi_t,d_{A_t}u_t) - \inner{\mu(u_t)}{\alpha_t}   
     \in\Om^1(\Sigma).      
$$
Moreover, for every $v\in\Lie(\Tilde{\Gg})$, 
$$
     \int_\Sigma h_vd\sigma_t
     = \int_\Sigma \sigma_t\wedge dh_v
     = \int_\Sigma \sigma_t\wedge\i(\pi_*v)\dvol_\Sigma
     = -\int_\Sigma \sigma_t(\pi_*v)\dvol_\Sigma.
$$
Hence
\begin{eqnarray*}
     \frac{d}{dt}\inner{\Tilde{\mu}(u_t,A_t)}{v}
&= &
     \int_\Sigma
     \bigl(
     \om(\xi_t,d_{A_t}u_t\circ v)
     + \inner{d\mu(u_t)\xi_t}{A_t(v)}
     \bigr)\dvol_\Sigma  \\
&&
     + \int_\Sigma 
       \bigl(
       \inner{A_t(v)}{d_{A_t}\alpha_t}
       + \inner{F_{A_t}}{\alpha_t(v)} 
       \bigr) \\
&= &
     \int_\Sigma
     \om(-du_t\circ v,\xi_t)\dvol_\Sigma  \\
&&
     + \int_\Sigma 
       \winner{(-d_{A_t}(A_t\circ v)-\i(v)F_{A_t})}{\alpha_t} \\
&= &
     \Om((-du_t\circ v,-\Ll_vA_t),(\xi_t,\alpha_t)).
\end{eqnarray*}
The last equality uses the formula
$
     \Ll_vA = d_A(A\circ v) + \i(v)F_A.
$
This proves the proposition.
\end{proof}

\begin{remark}
{\bf (i)}  
Consider the action of the gauge group $\Gg=\Gg(P)$ on $\Bb$.
The infinitesimal action of $\Om^0(\Sigma,\g_P)=\Lie(\Gg)$
is given by the vector fields
$$
    (u,A)\mapsto (L_u\eta,-d_A\eta)
$$
for $\eta\in\Om^0(\Sigma,\g_P)$.  
By Proposition~\ref{prop:Moment}, these vector fields
are Hamiltonian and a moment map
for the action is given by 
$$
     \Bb\to\Om^0(\Sigma,\g_P):
         (u,A)\mapsto *F_A+\mu(u).
$$
Hence the zero set of the moment map 
is the space of solutions of the second equation 
in~(\ref{eq:new}). 

\smallskip
\noindent{\bf (ii)}
If $M$ is K\"ahler then, 
under suitable regularity hypotheses,
the space 
$$
     \Xx = \left\{(u,A)\in\Bb\,|\,\bar\p_{J,A}(u)=0\right\}
$$
is a complex, and hence symplectic, submanifold 
of $\Bb$.  This submanifold is invariant under $\Gg$
and hence the space of gauge equivalence classes
of solutions of~(\ref{eq:new}) can be interpreted
as the symplectic quotient $\Xx\dslash\Gg$. 
In this case one can consider the action of the complexified 
group $\Gg^c$ on $\Xx$ and study the quotient
$\Xx^s/\Gg^c$, where $\Xx^s\subset\Xx$ is a suitable 
subspace of {\it stable pairs} $(u,A)$.  It turns out that,
as in the finite dimensional case, there is a natural 
correspondence
$$
     \Xx^s/\Gg^c\cong\Xx\dslash\Gg.
$$
This programme was carried out by Mundet in his recent
thesis~\cite{MUNDET}. 

\smallskip
\noindent{\bf (iii)}
The zero set of the moment map 
$\Tilde{\mu}:\Bb\to\Lie(\Tilde{\Gg})^*$
consists of all pairs $(u,A)\in\Bb$ that satisfy
$
     *F_A+\mu(u)=0
$
and
$$
     u^*\om-d\inner{\mu(u)}{A} 
     = \frac{\inner{[\om-\mu]}{B}}{\Vol(\Sigma)}
       \dvol_\Sigma,
$$
where the left hand side is understood as a 
$2$-form on $\Sigma$. 
However, the action of $\Tilde{\Gg}$ is not
compatible with the condition $\bar\p_{J,A}(u)=0$.
\end{remark}

In the case 
$$
     \Sigma=S^2
$$
with the standard metric and complex structure 
it is interesting to consider the group 
$\Hat{\Gg}\subset\Tilde{\Gg}$
of all automorphisms of $P$ that descend to 
isometries of $S^2$.  
There is an exact sequence 
$$
    1
    \TO\Gg
    \TO\Hat{\Gg}
    \TO\SO(3)
    \TO 1
$$
and the action of $\Hat{\Gg}$ on $\Bb$ preserves
the submanifold $\Xx$. 
The moment map 
$$
    \Hat{\mu}:\Bb\to\Lie(\Hat{\Gg})^*
$$ 
for this action is given by the restriction 
of $\Tilde{\mu}(u,A)$ to $\Lie(\Hat{\Gg})$.
Hence the zero set of $\Hat{\mu}$ consists
of all pairs $(u,A)\in\Bb$ that satisfy
$
     *F_A+\mu(u)=0
$
and
\begin{equation}\label{eq:donaldson1}
     \int_{S^2}
     h_\xi\bigl(
     u^*\om-d\inner{\mu(u)}{A}
     \bigr)
     = 0
\end{equation}
for every $\xi\in\so(3)$, where $h_\xi:S^2\to\R$
denotes the Hamiltonian function generating 
the infinitesimal action of $\xi$.  
Thus $h_\xi$ is the restriction of a linear
functional on $\R^3$ to $S^2$.  
It is interesting to consider all
solutions of~(\ref{eq:new}) and~(\ref{eq:donaldson1})
and divide by the action of the group $\Hat{\Gg}$.
This is the analogue of the quotient of the 
space of $J$-holomorphic spheres $v:S^2\to M$ 
by the reparametrization group $\PSL(2,\C)$.
Namely, $\PSL(2,\C)$ is the complexification of 
$\SO(3)$.  It acts on the Riemann sphere $\C\cup\{\infty\}$  
by fractional linear transformations of the form
$$
      \phi(z) = \frac{az+b}{cz+d},
$$
where $a,b,c,d$ are complex numbers such that $ad-bc=1$.
The subgroup of isometries is $\SO(3)\cong\SU(2)/\{\pm\1\}$.
The next proposition shows that,
instead of dividing the space of $J$-holomorphic 
curves from the Riemann sphere 
into a symplectic manifold by $\PSL(2,\C)$, 
one can consider $J$-holomorphic spheres that
satisfy~(\ref{eq:donaldson1}) (without the connection term)
and only divide by $\SO(3)$.
This is another analogue of the relation 
between the complex quotient $M^s/\G^c$
and the Marsden-Weinstein quotient $M\dslash\G$. 

\begin{proposition}\label{prop:balanced}
Let $\sigma\in\Om^2(S^2)$ such that
$$
      \int_{S^2}\sigma\ne 0.
$$
Then there exists a fractional linear transformation 
$\phi:S^2\to S^2$ such that
\begin{equation}\label{eq:balanced}
     \int_{S^2}h_\xi\phi^*\sigma = 0
\end{equation}
for every $\xi\in\so(3)$.
If $\phi_0\in\PSL(2,\C)$ and $\phi_1\in\PSL(2,\C)$ both 
satisfy~(\ref{eq:balanced}), and $\sigma$ is a volume form, 
then ${\phi_0}^{-1}\circ\phi_1\in\SO(3)$.
\end{proposition}

\begin{proof}
We identify the Riemann sphere $\C\cup\{\infty\}$ 
with the unit sphere $S^2\subset\R^3$ via stereographic
projection.  Explicitly, this diffeomorphism
is given by 
$$
     S^2\to\C\cup\{\infty\}:x\mapsto \frac{x_1+ix_2}{1-x_3}.
$$
Under this correspondence the quotient $\PSL(2,\C)/\SO(3)$
can be identified with the open unit ball $B^3\subset\R^3$
via the map 
$
    B^3\to\PSL(2,\C):
    \eta\mapsto\phi_\eta
$
that assigns to $\eta\in B^3$ the diffeomorphism
$\phi_\eta:S^2\to S^2$ given by 
$$
    \phi_\eta(x) 
    = \frac{\sqrt{1-|\eta|^2}}{1-\inner{x}{\eta}}
      \left(x-\inner{x}{\left|\eta\right|^{-1}\eta}
            \left|\eta\right|^{-1}\eta\right)
      + \frac{\inner{x}{\left|\eta\right|^{-1}\eta}-|\eta|}
             {1-\inner{x}{\eta}}\left|\eta\right|^{-1}\eta.
$$
If $\eta\in B^3$ converges to $\zeta\in S^2$ then
${\phi_\eta}^{-1}(x)$ converges to $\zeta$, 
uniformly in compact subsets of $S^2\setminus\{-\zeta\}$.
Hence 
$$
    \lim_{\eta\to\zeta}\int_{S^2} 
    \inner{\xi}{\i}{\phi_\eta}^*\sigma
    = \lim_{\eta\to\zeta}\int_{S^2} 
       \inner{\xi}{\i\circ{\phi_\eta}^{-1}}\sigma
    = \inner{\xi}{\zeta}\int_{S^2}\sigma
$$
for every $\zeta\in S^2$ and every $\xi\in\R^3$, 
where $\i:S^2\to\R^3$ denotes the obvious inclusion
and $\inner{\cdot}{\cdot}$ denotes the standard inner
product on $\R^3$.  Hence 
$$
    \lim_{\eta\to\zeta}\int_{S^2} 
    \i{\phi_\eta}^*\sigma
    = \zeta\int_{S^2}\sigma
$$
for every $\zeta\in S^2$, and the 
convergence is uniform in $\zeta$.  
It follows from a standard argument in degree theory~\cite{MILNOR}
that there exists an $\eta\in B^3$ such that
$$
     \int_{S^2} \i {\phi_\eta}^*\sigma = 0.
$$
This identity is equivalent to~(\ref{eq:balanced}).

Now suppose that $\sigma\in\Om^2(S^2)$ is a volume form
and $\phi\in\PSL(2,\C)$ such that
\begin{equation}\label{eq:phi-id}
    \int_{S^2} 
    \i\phi^*\sigma
    = \int_{S^2} 
      \i\sigma
    = 0.
\end{equation}
Then the same argument as in Remark~\ref{rmk:quotient}~(iii)
shows that $\phi\in\SO(3)$.  
Namely, there exists an $\eta\in B^3$ and a matrix
$\Psi\in\SO(3)$ such that 
$$
     \phi(x)=\phi_\eta(\Psi x)
$$
for every $x\in S^2$. 
Hence, for every $\xi\in\R^3$, 
$$
     \int_{S^2}\inner{\xi}{\i}{\phi_\eta}^*\sigma
     = \int_{S^2}\inner{\Psi^{-1}\xi}{\i}\phi^*\sigma
     = 0.
$$
Denote 
$
     \lambda(t)=\left|\eta\right|^{-1}\tanh(\left|\eta\right|t),
$
choose $T>0$ such that $\lambda(T)=1$, 
and consider the flow
$$
    \phi_t(x):=\phi_{\lambda(t)\eta}(x)
$$
for $0\le t\le T$.  It satisfies $\phi_0=\id$, $\phi_T=\phi_\eta$,
and, since $\dot\lambda=1-\left|\eta\right|^2\lambda^2$,
$$
    \frac{d}{dt}{\phi_t}^{-1}(x) 
    = \eta - \inner{\eta}{{\phi_t}^{-1}(x)}{\phi_t}^{-1}(x).
$$
Hence 
$$
    \frac{d}{dt}\int_{S^2} 
    \inner{\eta}{\i}{\phi_t}^*\sigma
= 
    \frac{d}{dt}\int_{S^2} 
    \inner{\eta}{\i\circ{\phi_t}^{-1}}\sigma  
= 
    \int_{S^2} 
        \bigl(\left|\eta\right|^2 
    - \inner{\eta}{\i\circ{\phi_t}^{-1}}^2
    \bigr)
    \sigma.
$$
The last expression is nonegative.
Integrating from $t=0$ to $t=T$, we obtain from~(\ref{eq:phi-id})
that it is equal to zero for all $t$. 
It follows that $\eta=0$, hence $\phi_\eta=\id$, 
and hence $\phi=\Psi\in\SO(3)$ as claimed.  
This proves the proposition. 
\end{proof}


\subsection{Hamiltonian perturbations}

We shall consider the following 
perturbations of~(\ref{eq:new}).
Let $\Cinf_\G(M)$ denote the space 
of $\G$-invariant smooth functions on $M$
and $\Vect_\G(M,\om)$ the space 
of $\G$-invariant Hamiltonian vector fields.
Choose a $\G$-invariant horizontal $1$-form 
$
     \sigma\in\Om^1(P,\Cinf_\G(M)).
$
One can think of $\sigma$ either as a $1$-form
on $\Sigma$ with values in $\Cinf_\G(M)$ 
or as a $1$-form 
on $P\times M$ that is invariant under the 
separate action of $\G$ on both $P$ and $M$
and that vanishes on all vectors of the form
$(p\xi,w)\in T_pP\times T_xM$, where $\xi\in\g$. 
Consider the $1$-form
$$
     X_\sigma:TP\to\Vect_\G(M,\om),\qquad
     \i(X_{\sigma_p(v)})\om = d(\sigma_p(v)).
$$
For every equivariant function $u:P\to M$
the $1$-form $X_\sigma(u)\in\Om^1(P,u^*TM)$,
defined by 
$
     T_pP\to T_{u(p)}M:v\mapsto X_{\sigma_p(v)}(u(p)),
$
is equivariant and horizontal.
Hence it descends to a $1$-form 
on $\Sigma$ with values in $u^*TM/\G$  
which will still be denoted by $X_\sigma(u)$.
The perturbed equations have the form
\begin{equation}\label{eq:new-ham}
     \bar\p_{J,A}(u) + (X_\sigma(u))^{0,1} = 0,\qquad
     *F_A + \mu(u) = 0.
\end{equation}
The space of solutions of~(\ref{eq:new-ham})
is invariant under the action
of the gauge group.

\begin{remark}\label{rmk:perturbation}
{\bf (i)}  
In local holomorphic coordinates~(\ref{eq:new-ham}) has the form
\begin{equation}\label{eq:new-ham-loc}
\begin{array}{rcl}
     \p_su+X_\Phi(u)+X_F(u)
     + J(\p_tu+X_\Psi(u)+X_G(u)) &= &0,  \\
     \p_s\Psi-\p_t\Phi+[\Phi,\Psi]
     +\lambda^2\mu(u) &= & 0,
\end{array}
\end{equation}
where $u:U\to M$, $\Phi,\Psi:U\to\g$, 
and $F,G:U\to\Cinf_\G(M)$. 
The local coordinate representatives
of $A$ and $\sigma$ are $A=\Phi\,ds+\Psi\,dt$ 
and $\sigma=F\,ds+G\,dt$. 

\smallskip
\noindent{\bf (ii)}
Let $\sigma\in\Om^1(P\times M)$ be as above.
Then $\sigma$ descends to a $1$-form on 
$P\times_\G M$ and $\om-d\inner{\mu}{A}-d\sigma$
is a connection $2$-form as in~\cite[Chapter~6]{MS2}.
The covariant derivative of a function $u:P\to M$
with respect to this connection is given by 
$$
     d_{A,\sigma}u
     = d_Au + X_\sigma(u)
     \in\Om^1(\Sigma,u^*TM/\G).
$$
The first equation in~(\ref{eq:new-ham})
can now be written in the form
$
     \bar\p_{J,A,\sigma}(u) = 0,
$
where $\bar\p_{J,A,\sigma}(u)$ is the $J$-antilinear
part of the $1$-form 
$d_{A,\sigma}u$.

\smallskip
\noindent{\bf (iii)}
The energy identity~(\ref{eq:energy})
continues to hold with $\bar\p_{J,A}(u)$
replaced by $\bar\p_{J,A,\sigma}(u)$
and $d_Au$ replaced by $d_{A,\sigma}u$
(in the definition of $E(u,A)$). 
\end{remark}


\subsection{Moduli spaces}

Fix an integer $k\ge 0$, a compact Riemann surface
$(\Sigma,J_\Sigma,\dvol_\Sigma)$, 
and an equivariant homology class 
$B\in H_2^\G(M;\Z)$ such that
\begin{equation}\label{eq:positive}
     \inner{[\om-\mu]}{B} > 0.
\end{equation}
By Proposition~\ref{prop:energy}
and Remark~\ref{rmk:perturbation}, 
this condition is necessary for the existence of solutions
of~(\ref{eq:new}) or~(\ref{eq:new-ham}). 
Consider the space
$$
     \Tilde{\Mm}_{B,\Sigma,k}
         = \Tilde{\Mm}_{B,\Sigma,k}(M,\mu;J,\sigma)
$$
of all tuples $(u,A,p_1,\dots,p_k)$ where $u$ and $A$ 
satisfy~(\ref{eq:new-ham}), $u$ represents the class $B$, 
and $p_1,\dots,p_k$ are points in $P$ with distinct
base points $\pi(p_i)\in\Sigma$. 
The gauge group $\Gg(P)$ acts on this space by
$$
     g^*(u,A,p_1,\dots,p_k)
         = (g^{-1}u,g^*A,p_1g(p_1)^{-1},\dots,p_kg(p_k)^{-1})
$$
and the quotient will be denoted by 
$$
     \Mm_{B,\Sigma,k}
     = \Mm_{B,\Sigma,k}(M,\mu;J,\sigma)
     = \Tilde{\Mm}_{B,\Sigma,k}(M,\mu;J,\sigma)/\Gg(P).
$$
If $k=0$ we write $\Mm_{B,\Sigma}=\Mm_{B,\Sigma,0}$.
Our goal is to use these moduli spaces to define invariants
of $(M,\om,\mu)$.   Note that the symplectic form 
enters in the definition of the moduli spaces only
indirectly through the compatibility condition on the 
almost complex structure $J$.


\subsection{Fredholm theory}

Let $\Bb$ be as in Section~\ref{subsec:red} 
and consider the vector bundle $\Ee\to\Bb$
whose fibre over a pair $(u,A)\in\Bb$ is given by
$
     \Ee_{u,A}
         = \Ee_u
         = \Om^{0,1}(\Sigma,u^*TM/\G) \oplus \Om^0(\Sigma,\g_P).
$
The gauge group $\Gg$ acts on $\Ee$, the projection 
$\Ee\to\Bb$ is $\Gg$-equivariant, 
and the almost complex structure $J$ and the perturbation $\sigma$
determine a $\Gg$-equivariant section 
\begin{equation}\label{eq:section}
     \Bb\to\Ee:
         (u,A)\mapsto (\bar\p_{J,A,\sigma}(u),*F_A+\mu(u)).
\end{equation}
Evidently, the zero set of this section
is the space $\Tilde{\Mm}_{B,\Sigma}$.
The vertical differential of the section~(\ref{eq:section})
at a zero $(u,A)$ gives rise to an 
operator
$$
     \Dd_{u,A}:
         \begin{array}{c}
                 \Cinf(\Sigma,u^*TM/\G) \\ 
                 \oplus  \\
                 \Om^1(\Sigma,\g_P)
     \end{array}
         \to 
         \begin{array}{c}
                 \Om^{0,1}(\Sigma,u^*TM/\G) \\
             \oplus  \\
                 \Om^0(\Sigma,\g_P)  \\
             \oplus  \\
                 \Om^0(\Sigma,\g_P)
         \end{array}
$$    
given by 
\begin{equation}\label{eq:DuA}
     \Dd_{u,A}
     \left(\begin{array}{c}\xi \\ \alpha\end{array}\right)
     = \left(
       \begin{array}{c}
         D\bar\p_{J,A,\sigma}(u)\xi + (L_u\alpha)^{0,1} \\
         {L_u}^*\xi - {d_A}^*\alpha  \\
         d\mu(u)\xi + *d_A\alpha
       \end{array}
       \right).
\end{equation}
Here 
$
     D\bar\p_{J,A,\sigma}(u):
         \Cinf(\Sigma,u^*TM/\G)
         \to \Om^{0,1}(\Sigma,u^*TM/\G)
$
is the Cauchy-Rie\-mann operator obtained by differentiating
the first equation in~(\ref{eq:new-ham}).  In explicit terms
this operator is given by 
$$
     D\bar\p_{J,A,\sigma}(u)\xi
     = (\Nabla{A,\sigma}\xi)^{0,1} 
       - \frac12J(\Nabla{\xi}J)\p_{J,A,\sigma}(u)
$$
for $\xi\in\Cinf(\Sigma,u^*TM/\G)$,
where $\nabla$ denotes the Levi-Civita connection
of the metric $\om(\cdot,J\cdot)$ on $M$ and
$\Nabla{A,\sigma}\xi\in\Om^1(\Sigma,u^*TM/\G)$ 
is given by
\begin{equation}\label{eq:nablaAsigma}
    \Nabla{A,\sigma}\xi 
    = \nabla\xi + \Nabla{\xi}X_A + \Nabla{\xi}X_\sigma.
\end{equation}

\begin{remark}\label{rmk:slice}
A tangent vector $(\xi,\alpha)\in T_{(u,A)}\Bb$
is $L^2$-orthogonal to the gauge orbit
of $(u,A)$ if and only if
\begin{equation}\label{eq:slice}
    {L_u}^*\xi - {d_A}^*\alpha = 0.
\end{equation}
This is the local slice condition and the tangent 
space of the quotient $\Bb/\Gg$ at $[u,A]$ can 
be identified with the space of solutions 
of~(\ref{eq:slice}).  Note also that the left hand side
of~(\ref{eq:slice}) agrees with the second coordinate
of $\Dd_{u,A}(\xi,\alpha)$ in~(\ref{eq:DuA}).
\end{remark}

The operator $\Dd_{u,A}$ (between suitable Sobolev 
completions) is a compact perturbation
of the direct sum of the first order operators
$D\bar\p_{J,A,\sigma}(u)$ and 
$$
     \Om^1(\Sigma,\g_P)
     \to\Om^0(\Sigma,\g_P)\oplus\Om^0(\Sigma,\g_P):
     \alpha\mapsto(-{d_A}^*\alpha,*d_A\alpha) 
$$
Hence it is Fredholm and, by the Riemann-Roch theorem,
\begin{equation}\label{eq:index}
     \INDEX\Dd_{u,A}
         = (2-2g)(n-\dim\G) + 2\inner{c_1^\G}{[u]}.
\end{equation}
Here $g$ is the genus of $\Sigma$ and
$
     c_1^\G = c_1^\G(TM,J)\in H^2_\G(M;\Z)
$
denotes the equivariant first Chern class of the 
tangent bundle.  It is defined as the first Chern class 
of the vector bundle $TM\times_\G\EG\to M\times_\G\EG$
with the complex structure given by $J\in\Jj(M,\om,\mu)$. 
If $\Dd_{u,A}$ is surjective for every 
$(u,A)\in\Tilde{\Mm}_{B,\Sigma}$
then it follows from the implicit function theorem
(in an infinite dimensional setting) that 
$\Mm_{B,\Sigma,k}$ is a smooth manifold
of dimension
\begin{equation}\label{eq:dim}
     \dim\,\Mm_{B,\Sigma,k}
         = (2-2g)(n-\dim\G) + 2\inner{c_1^\G}{B}
           + k(2+\dim\G).
\end{equation}
To obtain smooth moduli spaces it remains to prove that,
for a suitable perturbation, 
the Fredholm operator $\Dd_{u,A}$ is surjective 
for every solution $(u,A)$ of~(\ref{eq:new-ham}).
This means that the section~(\ref{eq:section})
of the vector bundle $\Ee\to\Bb$
is transverse to the zero section. 
In some cases transversality can be expected 
to hold for a generic volume form on $\Sigma$
or a generic almost complex structure on $M$. 
In other cases more general perturbations 
of the equations may be required. 
In~\cite{MUNDET} Mundet established transversality 
for the socalled {\it simple} (i.e. not multiply covered)
solutions in the case $\G=S^1$, by 
choosing a generic almost complex 
structure $J\in\Jj(M,\om,\mu)$.
Alternatively, denote by 
$\Mm_{B,\Sigma}^*\subset\Mm_{B,\Sigma}$
the subset of all those 
gauge equivalence classes
of solutions of~(\ref{eq:new-ham})
that satisfy 
\begin{equation}\label{eq:free}
     g\cdot u(p)=u(p)\qquad\IMP\qquad g=\1
\end{equation}
for every $p$ in a dense open subset of $P$.
The next proposition asserts that,
for a generic perturbation $\sigma$,
this subset is a manifold 
of the predicted dimension.
The proof will appear elsewhere.

\begin{proposition}\label{prop:trans}
Let $\Ss=\Cinf(\Sigma,\Cinf_\G(M))$ and denote
by $\Sreg\subset\Ss$ the subset of all perturbations
$\sigma\in\Ss$ such that the operator $\Dd_{u,A}$
is surjective for every solution $(u,A)$ of~(\ref{eq:new-ham})
that satisfies~(\ref{eq:free}).  
Then $\Sreg$ is a countable intersection of
dense open subsets of $\Ss$. 
\end{proposition}


\subsection{Compactness}

The moduli space $\Mm_{B,\Sigma}$ will in general not be 
compact.  The energy identity~(\ref{eq:energy}) 
asserts that 
$$
     E(u,A) = \inner{[\om-\mu]}{B}
$$
for every pair $(u,A)\in\Tilde{\Mm}_{B,\Sigma}$.
Hence the $L^2$-norms of $d_{A,\sigma}u$ 
and $\mu(u)$ are uniformly bounded. 
As in the case of $J$-holomorphic curves
and anti-self-dual Yang-Mills instantons,
this is a Sobolev borderline case.  
Combining the techniques 
for $J$-holomorphic curves (Gromov compactness~\cite{GRO}) 
with those for connections (Uhlenbeck compactness~\cite{UHL})
one can show that, for every sequence 
$(u^\nu,A^\nu)\in\Tilde{\Mm}_{B,\Sigma}$ that satisfies a 
uniform $L^p$ bound of the form
\begin{equation}\label{eq:Lp}
     \sup_\nu\int_\Sigma \bigl(
     |d_{A^\nu,\sigma}u^\nu|^p+|\mu(u^\nu)|^p
     \bigr)\dvol_\Sigma
     <\infty
\end{equation}
for some constant $p>2$, there exists a sequence
of gauge transformations $g^\nu\in\Gg(P)$ such that
$((g^\nu)^{-1}u^\nu,(g^\nu)^*A^\nu)$
has a $\Cinf$-convergent subsequence.
However, the energy identity only guarantees~(\ref{eq:Lp})
for $p=2$ and, in general, this does not suffice 
to prove compactness of the quotient space $\Mm_{B,\Sigma}$.
  
If~(\ref{eq:Lp}) does not hold then there
must be a sequence of points $p^\nu\in P$ 
such that either $|d_{A^\nu,\sigma}u^\nu(p^\nu)|$
or $|\mu(u^\nu(p^\nu))|$ diverges to infinity.
If the sequence $\mu\circ u^\nu$ is uniformly bounded
one can use the standard rescaling argument 
in Gromov compactness (cf.~\cite[Section~4.3]{MS1})
to prove that, for some sequence 
of maps 
$$
     \phi^\nu:\{z\in\C\,|\,|z|<1/\eps^\nu\}\to P
$$
with holomorphic projections $\pi\circ\phi^\nu$,
the sequence $u^\nu\circ\phi^\nu$ converges to
a nonconstant $J$-holomorphic curve $v:\C\to M$
that has finite energy.  The removable singularity
theorem for $J$-holomorphic curves 
(cf.~\cite[Theorem~4.2.1]{MS1}) then asserts 
that $v$ extends to a nonconstant 
$J$-holomorphic $2$-sphere in $M$. 
Any such $2$-sphere must be topologically 
nontrivial since 
$$
     E(v) = \int_{S^2}v^*\om = \int_{S^2} |dv|^2 > 0.
$$
Thus, if there are no $J$-holomorphic spheres in $M$
(for example if $\pi_2(M)=0$),
the only obstruction to compactness 
is the divergence of $\mu\circ u^\nu$. 
Now there are some interesting cases 
where the manifold $M$ is noncompact
but all solutions of~(\ref{eq:new}) satisfy 
a uniform bound on $u$. 

\begin{proposition}\label{prop:compact}
Assume the following.
\begin{description}
\item[(i)]
$(M,\om,J)$ is a Hermitian vector space.
\item[(ii)]
The group $\G$ acts on $M$ by unitary automorphisms.
\item[(iii)]
The moment map $\mu:M\to\g$ is proper. 
\end{description}
Then there exists a constant $c>0$ such that
every solution $(u,A)$ of~(\ref{eq:new})
(over any compact Riemann surface $\Sigma$) 
satisfies
\begin{equation}\label{eq:compact}
     \left\|u\right\|_{L^\infty}\le c.
\end{equation}
In particular, the moduli space 
$\Mm_{B,\Sigma}(M,\mu;J,\sigma)$ is compact
for every compact Riemann surface $(\Sigma,J_\Sigma,\dvol_\Sigma)$, 
every equivariant homology class $B\in H^\G_2(M;\Z)$,
and every compactly supported perturbation $\sigma$. 
\end{proposition}

\begin{proof}
Write $V:=M$, denote by 
$\inner{\cdot}{\cdot}=\om(\cdot,J\cdot)$
the real inner product on $V$,
by $\U(V)$ the group
of unitary automorphisms of $V$,
and by $\u(V)$ its Lie algebra.
By assumption, the action of $\G$ on $V$
is given by a homomorphism $\rho:\G\to\U(V)$
and we shall denote by $\dot\rho:\g\to\u(V)$
the corresponding Lie algebra homomorphism. 
We prove that there exists a central element
$\tau\in\g$ such that
\begin{equation}\label{eq:SW}
    \inner{x}{\dot\rho(\mu(x))Jx} 
    = 2\inner{\mu(x)}{\mu(x)-\tau}
\end{equation}
for $x\in V$. To see this,
suppose without loss of generality 
that $M=\C^n$ with its standard Hermitian 
structure, and consider the inner product 
$$
    \inner{A}{B} = \tr(A^*B)
$$
on the Lie algebra $\u(n)$ 
of skew-symmetric matrices.
The moment map is given by 
$$
    \mu(z) = \pi\left(-\frac{i}{2}zz^*\right) + \tau
$$
for $z\in\C^n$, where $\tau\in\g$ is a central element
and $\pi:\u(n)\to\g$ denotes the adjoint of the 
Lie algebra homomorphism $\dot\rho:\g\to\u(n)$.
Hence
$$
    \inner{z}{\dot\rho(\mu(z))iz}
    = \tr\left(\dot\rho(\mu(z))izz^*\right) 
    = \inner{\mu(z)}{\pi(-izz^*)} 
    = 2\inner{\mu(z)}{\mu(z)-\tau}.
$$
This proves~(\ref{eq:SW}).  

Now fix a Riemann surface $(\Sigma,J_\Sigma,\dvol_\Sigma)$
and suppose that $(u,A)$ is a solution of~(\ref{eq:new}).
Consider the equation~(\ref{eq:new-loc}) in local 
holomorphic coordinates, where the metric has
the form $\lambda^2(ds^2+dt^2)$.  In our situation
$$
    X_\xi(u)=\dot\rho(\xi)u
$$
and we abbreviate
$
    \Nabla{s}u := \p_su+\dot\rho(\Phi)u
$
and
$
    \Nabla{t}u := \p_tu+\dot\rho(\Psi)u.
$
Since
$$
    \Nabla{s}\Nabla{t}u-\Nabla{t}\Nabla{s}u 
    = \dot\rho(\p_s\Psi-\p_t\Phi+[\Phi,\Psi])u
$$
and, by~(\ref{eq:new-loc}), 
$$
    \Nabla{s}u+J\Nabla{t}u=0,\qquad
    \p_s\Psi-\p_t\Phi+[\Phi,\Psi] + \lambda^2\mu(u)=0,
$$ 
we obtain
$$
    \Nabla{s}\Nabla{s}u+\Nabla{t}\Nabla{t}u
    = \lambda^2\dot\rho(\mu(u))Ju.
$$
Hence, with $\Delta={\p_s}^2+{\p_t}^2$, 
\begin{eqnarray*}
    \Delta |u|^2/2
&= &
    \p_s\inner{u}{\Nabla{s}u} + \p_t\inner{u}{\Nabla{t}u} \\
&= &
    \left|\Nabla{s}u\right|^2 + \left|\Nabla{t}u\right|^2 
    + \inner{u}{\Nabla{s}\Nabla{s}u+\Nabla{t}\Nabla{t}u} \\
&= &
    \left|\Nabla{s}u\right|^2 + \left|\Nabla{t}u\right|^2 
    + \lambda^2\inner{u}{\dot\rho(\mu(u))Ju} \\
&= &
    \left|\Nabla{s}u\right|^2 + \left|\Nabla{t}u\right|^2 
    + 2\lambda^2\inner{\mu(u)}{\mu(u)-\tau} \\
&\ge &
    2\lambda^2\left|\mu(u)\right|
    \left(\left|\mu(u)\right|-\left|\tau\right|\right)
\end{eqnarray*}
Now let $(s_0,t_0)$ be a point at which the function
$(s,t)\mapsto|u(s,t)|$ attains its maximum.
Since $\Sigma$ is compact, such a point exists 
in some coordinate chart, and we have
$\Delta|u|^2\le 0$ at $(s_0,t_0)$.  Hence 
$$
    \left|\mu(u(s_0,t_0))\right|
    \le \left|\tau\right|.
$$
Since $\mu$ is proper, there exists a constant 
$c>0$ such that
$$
    |\mu(x)|\le|\tau|
    \qquad\IMP\qquad
    |x|\le c.
$$
Hence $|u(s_0,t_0)|\le c$ and it follows that
$
    \sup_{p\in P}|u(p)| \le c
$
for every solution $(u,A)$ of~(\ref{eq:new}).
To prove the last assertion just note that the same 
estimate holds for solutions of the perturbed 
equation~(\ref{eq:new-ham}) whenever the 
support of the perturbation is contained
in the ball $\{|x|<c\}$. 
This proves the proposition.
\end{proof}

\begin{remark}\label{rmk:bubble}
{\bf (i)}   The proof of Proposition~\ref{prop:compact}
is reminiscent of the compactness proof for the 
Seiberg--Witten equations in 
Kronheimer--Mrowka~\cite{KM}.

\smallskip
\noindent{\bf (ii)}
In Proposition~\ref{prop:compact} the assumption 
that the moment map be proper is essential.
But one would expect that conditions~(i) and~(ii)
can be removed or be replaced by weaker 
assumptions.   

\smallskip
\noindent{\bf (iii)}
If $\pi_2(M)\ne 0$ then, in general, there may be 
$J$-holomorphic spheres in $M$. 
In this case the compactification of the 
moduli space $\Mm_{B,\Sigma}(M,\mu;J,\sigma)$ 
should include stable maps, as introduced
by Kontsevich~\cite{KONTSEVICH}. To see this think
of the solutions of the first equation in~(\ref{eq:new-ham}) 
as $\Tilde{J}_{A,\sigma}$-holomorphic
curves from $\Sigma$ to $\Tilde{M}$
(see page~\pageref{tildeJA} for the case $\sigma=0$).
In the stable maps that appear in the 
limit the main component will be a solution
of~(\ref{eq:new-ham}) and all other components 
will be $J$-holomorphic spheres in the fibres.

\smallskip
\noindent{\bf (iv)}
It is often interesting to allow the complex 
structure on $\Sigma$ to vary.  One then has to deal 
with a suitable compactification of Teichm\"uller 
space.  This again leads to Kontsevich's stable maps.

\smallskip
\noindent{\bf (v)}
Similar techniques as in the proof 
of Proposition~\ref{prop:compact}
can be used to prove the following unique continuation
theorem.  
{\it
Let $(u,A)$ be a solution of~(\ref{eq:new-ham}).
If the pair $(d_{A,\sigma}u,\mu\circ u)$ 
vanishes to infinite order
at some point $p\in P$ then $d_{A,\sigma}u\equiv0$ 
and $\mu(u)\equiv0$.
}
\end{remark}


\subsection{Invariants}

The moduli space $\Tilde{\Mm}_{B,\Sigma,k}(M,\mu;J,\sigma)$
carries a natural right action of 
$\G^k=\G\times\cdots\times\G$ on the 
$k$ marked points.  This action commutes with the 
action of the gauge group and hence descends
to an action on the quotient space 
$\Mm_{B,\Sigma,k} = \Mm_{B,\Sigma,k}(M,\mu;J,\sigma)$.
Now there is an evaluation map
$\ev=(\ev_1,\dots,\ev_k):\Mm_{B,\Sigma,k}\to M^k$,
given by 
$
     \ev_i([u,A,p_1,\dots,p_k]) = u(p_i)
$
and a projection $\pi:\Mm_{B,\Sigma,k}\to\Aa(P)/\Gg(P)$
given by 
$
     \pi([u,A,p_1,\dots,p_k])=[A].
$
The evaluation map is $\G^k$-equivariant 
and the projection $\pi$ is $\G^k$-invariant.
$$
\begin{array}{ccc}
    \Mm_{B,\Sigma,k} & \stackrel{\ev}{\TO} & M^k \\
  {\scriptstyle{^\pi}}\downarrow\;\;\, & & \\
   \Aa/\Gg & & 
\end{array}.
$$
One can use these maps to produce certain 
natural $\G^k$-equivariant cohomology classes on the 
moduli space $\Mm_{B,\Sigma,k}$.  Integrating these 
over the quotient $\Mm_{B,\Sigma,k}/\G^k$ gives rise 
to the invariants. 

To be more precise choose equivariant cohomology 
classes $\alpha_i\in H^*_\G(M)$ for $i=1,\dots,k$
and a cohomology class $\beta\in H^*(\Aa/\Gg)$ such that
\begin{equation}\label{eq:deg}
    \deg(\beta) + \sum_{i=1}^k\deg(\alpha_i)
    = \dim\Mm_{B,\Sigma,k} - k\dim\G.
\end{equation}
Then the pullback 
$
    \pi^*\beta\smile{\ev_1}^*\alpha_1\smile\cdots\smile{\ev_k}^*\alpha_k
$
is an equivariant cohomology class on $\Mm_{B,\Sigma,k}$. 
Let us pretend, for a moment, that $\Mm_{B,\Sigma,k}$ 
is a compact smooth manifold of the predicted dimension
and that $\G^k$ acts freely on this space.  
Then our equivariant cohomology class on $\Mm_{B,\Sigma,k}$
descends to a top dimensional
cohomology class on the quotient $\Mm_{B,\Sigma,k}/\G^k$
that we can evaluate on the fundamental cycle.
This gives rise to an integer
\begin{equation}\label{eq:invariant}
    \Phi_{B,\Sigma,k}^{M,\mu}(\beta,\alpha_1,\dots,\alpha_k)
    := \int_{\Mm_{B,\Sigma,k}/\G^k}
    \pi^*\beta\smile{\ev_1}^*\alpha_1\smile\cdots\smile{\ev_k}^*\alpha_k.
\end{equation}
In general, only the subspace $\Mm_{B,\Sigma,k}^*$ of all solutions
that satisfy~(\ref{eq:free}) for almost every $p\in P$ 
and for $p=p_i$ is a smooth manifold
for a generic perturbation $\sigma$ 
and carries a free action of $\G^k$.  
Even under the hypotheses
of Proposition~\ref{prop:compact} 
this space will not be compact. 
However, in many cases we expect that 
this space can be compactified by 
adding strata of strictly lower dimensions,
and that~(\ref{eq:invariant}) can be defined
by integrating differential forms whose 
pullbacks are supported in $\Mm_{B,\Sigma,k}^*$.  
Alternatively, one can consider intersection 
numbers of cycles in $M\times_\G\EG$.
This requires the choice of an equivariant 
function $\Mm_{B,\Sigma,k}^*\to\EG$ and the easiest
way to get such a function is by composition
of the evaluation map with an equivariant
function $\phi:M^*\to\EG$, where
$$
     M^*=\{x\in M\,|\,gx=x\IMP g=\1\}.
$$
Here we assume that $\EG$ has been replaced
by a suitable finite dimensional approximation. 
Now represent the Poincar\'e duals of $\alpha_i$ and $\beta$
by submanifolds $Y_i\subset M\times_\G\EG$ 
and $Z\subset\Aa/\Gg$. Integrating the differential 
form then corresponds to counting the solutions
$[u,A,p_1,\dots,p_k]\in\Mm_{B,\Sigma,k}^*$
that satisfy
\begin{equation}\label{eq:Yi}
    [u(p_i),\phi(u(p_i))]\in Y_i,\qquad
        [A]\in Z.
\end{equation}
Since $\phi$ is only defined on $M^*$ one
has to check that the reducible solutions 
of~(\ref{eq:new-ham}), if they exist, 
do not obstruct compactness.
With standard cobordism techniques,
similar to the ones used in the 
definition of the Donaldson invariants~\cite{DONALDSON2},
the Gromov--Witten invariants~\cite{MS1,RUAN2},
or the Seiberg--Witten invariants~\cite{SAL1},
one should then be able to prove 
that the invariants~(\ref{eq:invariant})
are independent of the choice of the perturbation $\sigma$
and the almost complex structure $J$ used to define them.
To work this out in detail requires 
a considerable amount of analysis
which will be carried out elsewhere.  
Some cases were treated by Mundet~\cite{MUNDET}.

\begin{remark}\label{rmk:stable-curve}
In the above discussion the complex structure on the 
Riemann surface $\Sigma$ is fixed.  Even in this case 
there is an interesting moduli space $\Mm_{\Sigma,k}$ 
of stable Riemann surfaces with $k$ marked points,
where one of the components of the stable surface 
is $\Sigma$ itself.  Correspondingly, one might wish
to extend the definition of the invariants to 
include, as a base for the bundle $P$, stable 
Riemann surfaces where the main component 
is $\Sigma$ and all other components are spheres.
With this modification in place 
there is a projection 
$$
     \Mm_{B,\Sigma,k}\to\Mm_{\Sigma,k}
$$
and one could consider pullbacks of cohomology classes 
from $\Mm_{\Sigma,k}$ to get further invariants.  
Similar observations apply to the case where the 
complex structure on $\Sigma$ is allowed to vary.
\end{remark}


\subsection{Adiabatic limits}\label{sec:adiabatic}

In her PhD thesis~\cite{GAIO} 
the second author studies the adiabatic 
limit $\eps\to 0$ in the equations
\begin{equation}\label{eq:rita}
    \bar\p_{J,A}(u) = 0,\qquad
    *F_A + \eps^{-2}\mu(u) = 0.
\end{equation}
For $\eps=0$ these equations degenerate into~(\ref{eq:jhol})
and the solutions of those equations correspond to 
$J$-holomorphic curves in the Marsden-Weinstein quotient
$M\dslash\G$ (see Section~\ref{sec:jhol}).  
Under suitable conditions on $M$ and for 
sufficiently small $\eps>0$, there should be a one-to-one 
correspondence between the solutions 
of~(\ref{eq:rita}) and those of~(\ref{eq:jhol}).  
The arguments are reminiscent of the proof 
of the Atiyah--Floer conjecture in~\cite{DOSAL3,SAL4}.
Gaio proves that regular solutions of~(\ref{eq:jhol}) 
give rise to solutions of~(\ref{eq:rita}) 
for $\eps$ sufficiently small,
and makes substantial progress towards establishing 
that, in many cases, all solutions of~(\ref{eq:rita})
can be obtained in this way.  
When completed, this work should
lead to a proof of the following conjecture,
at least in the case where the quotient $M\dslash\G$ 
is {\it semi-positive} (or {\it weakly monotone}
in the terminology of~\cite{HOSA,MS1}).

\begin{conjecture}\label{con:adiabatic}
Suppose that $\mu:M\to\g$ is proper,
that $0$ is a regular value of $\mu$,
that $\mu^{-1}(0)$ is nonempty, 
and that $\G$ acts freely on $\mu^{-1}(0)$.
Then, for $\bar B\in H_2(M\dslash\G;\Z)$
and $\alpha_1,\dots,\alpha_k\in H^*_\G(M;\Z)$,
\begin{equation}\label{eq:gromov}
     \Phi_{B,\Sigma,k}^{M,\mu}(1,\alpha_1,\dots,\alpha_k)
     = {\rm GW}_{\bar B,\Sigma,k}^{M\dslash\G}
       (\bar\alpha_1,\dots,\bar\alpha_k).
\end{equation}
Here $1\in H^0(\Aa/\Gg)$,
$B\in H^\G_2(M;\Z)$ is the image of 
$\bar B$ under the homomorphism 
$H_2(M\dslash\G;\Z)\to H^\G_2(M;\Z)$ 
induced by the inclusion $\mu^{-1}(0)\INTO M$,
and $\bar\alpha_i\in H^*(M\dslash\G;\Z)$ 
is the image of $\alpha_i$ under 
the homomorphism $H^*_\G(M;\Z)\to H^*(M\dslash\G;\Z)$ 
induced by the same inclusion.
\end{conjecture}

\begin{remark}\label{rmk:adiabatic}
{\bf (i)}
Kirwan~\cite{MFK} proved that the homomorphism
$
     H^*_\G(M;\Z)\to H^*(M\dslash\G;\Z)
$
is surjective.

\smallskip
\noindent{\bf (ii)}
Consider the case
$$
     M\dslash\G = \{\point\}.
$$
Then $n=\dim\G=\dim M/2$ and the only class
in the image of the homomorphism 
$H_2(M\dslash\G;\Z)\to H^\G_2(M;\Z)$ 
is $B=0$.  Moreover, the invariant 
$\Phi_{B,\Sigma,k}^{M,\mu}(1,\alpha_1,\dots,\alpha_k)$
can only be nonzero if $\dim\Mm_{B,\Sigma,k}=0$.
Hence assume
$$
     n=\dim\G,\qquad B=0,\qquad k=0.
$$
Then Conjecture~\ref{con:adiabatic} asserts that, 
if $0$ is a regular value of $\mu$ and
$\G$ acts freely on $\mu^{-1}(0)$, then
$$
     \Phi_{0,\Sigma,0}^{M,\mu} = 1.
$$
In this case the bundle $P$ is trivial
and, for any $\eps>0$, there are obvious solutions 
of~(\ref{eq:rita}) that satisfy $\mu\circ u=0$,
$d_Au=0$, and $F_A=0$.  They are all gauge 
equivalent.  The proof of Conjecture~\ref{con:adiabatic} 
would be to show that, for $\eps>0$ sufficiently small,
there is no other solution of~(\ref{eq:rita}).

\smallskip
\noindent{\bf (iii)}
Conjecture~\ref{con:adiabatic} does not allow
for the pullback of classes in $\Mm_{\Sigma,k}$
or for variations of the complex structure 
on $\Sigma$ (see Remark~\ref{rmk:stable-curve}).
But there should be analogous results for 
those cases.   

\smallskip
\noindent{\bf (iv)}
The examples in Sections~\ref{sec:vortex}
and~\ref{sec:bradlow} show that the 
invariants $\Phi_{B,\Sigma,k}^{M,\mu}$ can be nontrivial
in cases where the symplectic quotient $M\dslash\G$ 
is a point or the empty set 
(and $B$ does not descend to a homology 
class in the quotient). 

\smallskip
\noindent{\bf (v)}
If $\G$ does not act freely on $\mu^{-1}(0)$
then Conjecture~\ref{con:adiabatic} suggests that the 
solutions of~(\ref{eq:new-ham}) can be used
to define the Gromov--Witten invariants 
of symplectic orbifolds.
\end{remark}


\subsection{Wall crossing and localization}\label{sec:wall}

One should be able to use the formula~(\ref{eq:gromov}) to find 
relations between the Gromov--Witten invariants of the quotients
$M\dslash\G(\tau)$ for different values of~$\tau$.
Namely, choose a generic path 
$
     [0,1]\to\g:s\mapsto\tau_s
$
in the center of $\g$ and consider the cobordism 
$$
     \Ww_{B,\Sigma} = \bigcup_{0\le s\le 1}
     \{s\}\times\Mm_{B,\Sigma}(\mu-\tau_s)
$$
with boundary
$$
     \p\Ww_{B,\Sigma} 
     = \Mm_{B,\Sigma}(\mu-\tau_0)
       \cup \Mm_{B,\Sigma}(\mu-\tau_1).
$$
In some cases the critical parameters should 
be the singular values of the moment map 
(e.g. when $\G=S^1$).  
However, the examples in 
Sections~\ref{sec:bradlow} and~\ref{sec:grass} show
that the moduli space $\Mm_{B,\Sigma}(\mu-\tau)$
may also have singularities when $\tau$ is a regular value
of the moment map, and the effect of these on the
definition of the invariants remains yet to be fully
understood.  If the path $s\mapsto\tau_s$ passes 
through such critical parameters then the difference of the 
invariants for $\tau_0$ and $\tau_1$ 
should be computable in terms of the 
reducible solutions of~(\ref{eq:new-ham}).

\begin{remark}\label{rmk:wall}
{\bf (i)}
For the ordinary cohomology of symplectic quotients 
wall crossing formulae were discovered by 
Martin~\cite{MARTIN,MARTIN1,MARTIN2}. 
These should correspond to the present case
when $B=0$. 
In~\cite{MARTIN2} Martin developed techniques for 
computing the cohomology of symplectic
quotients via a reduction argument to the action 
of the maximal torus.  We expect that his ideas 
can be adapted to our situation and lead
to formulae for the computation of the invariants 
$\Phi_{B,\Sigma,k}$.

\smallskip
\noindent{\bf (ii)}
Guillemin and Sternberg~\cite{GS} showed 
that passing through a critical value of 
the moment map corresponds to blowing 
up and down.  
The resulting formulae should thus lead 
to an alternative proof of Ruan's 
results in~\cite{RUAN1}.

\smallskip
\noindent{\bf (iii)}
We expect that the wall crossing relations 
correspond, under suitable assumptions,
to the fixed point localization 
formulae of Kontsevich~\cite{KONTSEVICH} 
and Givental~\cite{GIVENTAL1}.
In~\cite{GIVENTAL2,GIVENTAL3} 
Givental used localization 
to compute Gromov--Witten invariants 
for many examples and, in particular, 
to prove the mirror conjecture for the 
quintic in $\C P^4$. 

\smallskip
\noindent{\bf (iv)}
It should be interesting 
to relate the wall-crossing formulae 
for the invariants $\Phi_{B,\Sigma,k}$ 
to the gluing formulae in contact homology
(cf. Eliashberg--Givental--Hofer~\cite{EGH}).
\end{remark}


\section{Floer homology}\label{sec:floer}


\subsection{Relative fixed points}

Let $(M,\om,\mu)$ be a symplectic manifold
with a Hamiltonian $\G$-action. 
Fix a time dependent Hamiltonian function
$
     \R\times M\to\R:(t,x)\mapsto H_t(x)
$
such that 
$
     H_t=H_{t+1}
$
and $H_t:M\to\R$ is $\G$-invariant for every $t$. 
Consider the Hamiltonian differential equation
\begin{equation}\label{eq:ham}
     \dot x(t) = X_{H_t}(x(t))
\end{equation}
and denote by $f:M\to M$ the time-$1$ map.
It is defined by $f(x(0))=x(1)$ for all solutions
of~(\ref{eq:ham}).  Note that 
$$
     \mu\circ f=\mu.
$$
A pair $(x_0,g_0)\in M\times\G$ is called 
is called a {\bf relative fixed point} of $f$ if
$$
     f(x_0) = g_0x_0.
$$
Equivalently, the unique solution 
$x:\R\to M$ of~(\ref{eq:ham}) with initial
condition $x(0)=x_0$ satisfies 
$x(t+1)=g_0x(t)$ for every $t\in\R$. 
Note that the set of relative fixed points
is invariant under the action of $\G$ on
$M\times\G$ by $(x_0,g_0)\mapsto(gx_0,gg_0g^{-1})$.
A relative fixed point $(x_0,g_0)$ is called
{\bf regular} if $gx_0=x_0$ implies $g=\1$.
It is called {\bf nondegenerate} if the linear map 
$
     df(x_0)-g_0:T_{x_0}M\to T_{g_0x_0}M
$
induces an isomorphism from the quotient
$\ker d\mu(x_0)/L_{x_0}\g_\tau$
to $\ker d\mu(g_0x_0)/L_{g_0x_0}\g_\tau$,
where $\tau=\mu(x_0)$. This means that  
\begin{equation}\label{eq:nondegenerate}
     d\mu(x_0)v=0,\quad
     df(x_0)v - g_0v\in\im L_{g_0x_0}\qquad\IMP\qquad
     v\in\im L_{x_0}
\end{equation}
for every $v\in T_{x_0}M$. 
Relative fixed points $(x_0,g_0)\in\mu^{-1}(0)\times\G$ 
appear as the critical points of an equivariant 
symplectic action functional.  


\subsection{Equivariant symplectic action}

Denote by $D\subset\C$ the closed 
unit disc and by 
$
     \Ll = \Ll(M\times\g)
$
the space of contractible loops in $M\times\g$.
The universal cover of this space consists 
of all equivalence classes of triples
$(x,\eta,v)$, where $x:\R/\Z\to M$, $\eta:\R/\Z\to\g$
and $v:D\to M$ satisfy $v(e^{2\pi it})=x(t)$. 
Two such triples $(x_1,\eta_1,v_1)$ and 
$(x_2,\eta_2,v_2)$ are equivalent iff
$x_1=x_2$, $\eta_1=\eta_2$, and $v_1$ is homotopic
to $v_2$ with fixed boundary.  The space
of equivalence classes will be denoted by 
$$
     \Tilde{\Ll} = \Tilde{\Ll}(M\times\g).  
$$
This space carries an action of the group
$
     \Tilde{\Gg}=\Map(D,\G)
$
by 
$$
     g^*[x,\eta,v] = [g^{-1}x,g^{-1}\p_tg + g^{-1}\eta g,g^{-1}v],
$$
where $\p_tg = \p/\p t\,g(e^{2\pi it})$.
There is a $\Tilde{\Gg}$-invariant 
action functional
$$
    \Aa_{\mu,H}:\Tilde{\Ll}(M\times\g)\to\R
$$
given by 
$$
    \Aa_{\mu,H}(x,\eta,v)
    = - \int_Dv^*\om 
      + \int_0^1\bigl(
        \inner{\mu(x(t))}{\eta(t)} - H_t(x(t))
        \bigr)\,dt.
$$
A $1$-periodic family of almost complex structures 
$J_t\in\Jj(M,\om,\mu)$ determines
an $L^2$-inner product on the tangent space
$$
     T_{(x,\eta)}\Ll
     = \Cinf(S^1,x^*TM)\times\Cinf(S^1,\g),
$$
and the gradient of $\Aa_{\mu,H}$ with respect
to this inner product is given by 
\begin{equation}\label{eq:grad}
     \grad\Aa_{\mu,H}(x,\eta)
     = \left(\begin{array}{c}
       J_t\left(\dot x+X_\eta(x)-X_{H_t}(x)\right)  \\
           \mu(x)
       \end{array}
       \right).
\end{equation}
Hence the critical points of $\Aa_{\mu,H}$ are the
loops $(x,\eta):\R/\Z\to M\times\g$ that satisfy
\begin{equation}\label{eq:relative}
     \dot x + X_\eta(x) = X_{H_t}(x),\qquad
     \mu(x) = 0.
\end{equation}
Let us denote by $\Tilde{\Per}(\mu,H)$
the set of solutions of~(\ref{eq:relative}).
The loop group
$$
     L\G = \Map(S^1,\G)
$$
acts on this space and the quotient will be denoted by 
$$
     \Per(\mu,H) = \Tilde{\Per}(\mu,H)/L\G.
$$
This quotient space can be naturally
identified with the set of $\G$-orbits 
of relative fixed points 
of $f$ in $\mu^{-1}(0)\times\G$.  Moreover, 
a relative fixed point $(x_0,g_0)$ is nondegenerate
if and only if the corresponding critical
point of $\Aa_{\mu,H}$ is nondegenerate.
The proof of this observation is a precise analogue 
of the proof of Proposition~4.4 in~\cite{DOSAL2}.

\begin{remark}
A closer look at the equivariant symplectic 
action should reveal interesting
relations to the geometry of
the loop group (cf.~\cite{PS,DAVIES}).
\end{remark}


\subsection{Floer homology}

One can construct Floer homology
groups $\HF^*(M,\om,\mu;J,H)$,
as in the standard case,  
by considering the gradient 
flow lines of the action functional $\Aa_{\mu,H}$ 
with respect to the $L^2$-metric determined by $J_t$.
The formula~(\ref{eq:grad}) shows that the gradient
flow lines are pairs $(u,\Psi)$ where 
$u:\R^2\to M$ and $\Psi:\R^2\to\g$ satisfy
\begin{equation}\label{eq:floer}
     \p_su+J_t(\p_tu+X_\Psi(u)-X_{H_t}(u)) = 0,\qquad
     \p_s\Psi + \mu(u) = 0,
\end{equation}
and
\begin{equation}\label{eq:per}
     u(s,t+1) = u(s,t),\qquad 
     \Psi(s,t+1) = \Psi(s,t). 
\end{equation}
The energy of such a flow line is defined by 
$$
    E(u,\Psi)
    = \int_0^1\int_{-\infty}^\infty
      \bigl(
      \left|\p_su\right|^2 + \left|\p_s\Psi\right|^2
      \bigr)
      \,dsdt.
$$
If this energy is finite, and the critical points
of $\Aa_{\mu,H}$ are all nondegenerate, then one can
show with standard techniques in gauge theory
that the limits
\begin{equation}\label{eq:limit}
     x^\pm(t) = \lim_{s\to\pm\infty}u(s,t),\qquad
     \eta^\pm(t) = \lim_{s\to\pm\infty}\Psi(s,t) 
\end{equation}
exist and are critical points of $\Aa_{\mu,H}$. 
The strategy would now be to proceed as in the standard
case and define a chain complex generated by 
the critical points of $\Aa_{\mu,H}$ and 
define a boundary operator by counting 
the solutions of~(\ref{eq:floer})
and~(\ref{eq:per}) with given limits~(\ref{eq:limit})
in the case where the Floer relative Morse index is $1$. 
To carry this out in detail one has to deal
with the usual transversality and compactness 
questions. Additional difficulties arise from the presence 
of nontrivial isotropy subgroups of critical
points of $\Aa_{\mu,H}$ and this would require
an equivariant version of Floer homology 
(cf. Viterbo~\cite{VITERBO}). 
The resulting Floer homology theory is related
to the solutions of~(\ref{eq:new}) in the same way
as instanton Floer homology~\cite{FLOER1}
is related to the Donaldson invariants,
symplectic Floer homology~\cite{FLOER2} 
is related to the Gromov--Witten invariants,
and Seiberg--Witten Floer homology is related to
the Seiberg--Witten invariants.
To see this compare~(\ref{eq:floer}) 
with~(\ref{eq:new-ham-loc}).
In particular, one should get relative invariants, 
for Riemann surfaces with cylindrical ends, 
with values in the Floer homology groups
(cf.~\cite{PSS} for the standard case).  


\subsection{The Arnold conjecture for regular quotients}

Suppose that $\mu$ is proper, 
$0$ is a regular value of $\mu$, 
and $\G$ acts freely on $\mu^{-1}(0)$.  
Then one would hope to obtain transversality
for the solutions of~(\ref{eq:floer}) by choosing
a generic $\G$-invariant Hamiltonian $H$.
At first glance one might not expect to get 
anything new, because the critical points 
are the periodic solutions 
of a Hamiltonian system in $M\dslash\G$ 
and one could get an equivalent theory 
from Floer homology in the reduced space.
However, the compactness result of
Proposition~\ref{prop:compact} suggests 
that in many cases the present approach might
be simpler than the standard theory,
and lead to a proof of the Arnold conjecture
over the integers.\footnote
{
A quite different approach to Floer 
homology over the integers for general 
symplectic manifolds has recently 
been proposed by Fukaya~\cite{FUKAYA}.
Other approaches to Floer homology
for general symplectic manifolds
(cf. Fukaya-Ono~\cite{FO}, Liu-Tian~\cite{LT})
have so far only been established 
over the rationals.
}
The key point is that the presence of 
holomorphic spheres with negative Chern number
in the quotient $M\dslash\G$ leads to
complications in the standard theory,
but not in our approach, provided that 
they do not lift to holomorphic spheres in $M$. 


\subsection{Relation with Morse theory}

If $0$ is not a regular value of $\mu$,
or $\G$ does not act freely on $\mu^{-1}(0)$,
then the Floer homology theory outlined above
should lead to new existence theorems for
relative fixed points of $\G$-equivariant
Hamiltonian symplectomorphisms.
In the standard theory one can, in many cases, 
identify the Floer homology groups with Morse 
homology by considering the case where $H$
is independent of $t$ (and sufficiently small). 
The analogue of this argument in the present 
case leads to equivariant Morse homology
on $M\times\g$ for the function
$$
     M\times\g\to\R:(x,\eta)\mapsto \inner{\mu(x)}{\eta}-H(x),
$$
where $H:M\to\R$ is $\G$-invariant. The 
critical points of this function satisfy
$$
     \nabla H(x) = JX_\eta(x),\qquad \mu(x)=0,
$$
and so correspond to critical points of the 
induced function $\bar H:M\dslash\G\to\R$
whenever the quotient is smooth. 
The gradient flow equations have the form
\begin{equation}\label{eq:gradflow}
     \dot u + JX_\Psi(u)-\nabla H(u) = 0,\qquad
     \dot\Psi + \mu(u) = 0.
\end{equation}
They are equivalent to~(\ref{eq:floer})
whenever $H$, $J$, $u$, and $\Psi$ are independent of~$t$.


\subsection{Equivariant symplectomorphisms}

One might wish to define the Floer homology
groups of general equi\-variant symplectomorphisms,
not just Hamiltonian ones. For this theory one would 
consider symplectomorphisms $f:M\to M$ 
that satisfy
\begin{equation}\label{eq:mu-phi}
     f(gx) = \rho(g)x,\qquad
     \mu(f(x)) = \dot\rho(\mu(x)),
\end{equation}
for all $x\in M$ and some isomorphism $\rho:\G\to\G$.  
Here $\dot\rho:\g\to\g$ denotes the corresponding 
Lie algebra isomorphism. 
The Hamiltonian perturbation
$H_t\in\Cinf_\G(M)$ and the almost complex structures
$J_t\in\Jj(M,\om,\mu)$ should satisfy the periodicity
condition
$$
     H_t = H_{t+1}\circ f,\qquad
     J_t = f^*J_{t+1},
$$
and~(\ref{eq:per}) should be replaced by
\begin{equation}\label{eq:per-phi}
     u(s,t+1) = f(u(s,t)),\qquad 
     \Psi(s,t+1) = \dot\rho(\Psi(s,t)). 
\end{equation}
The resulting solutions of~(\ref{eq:floer}) 
and~(\ref{eq:per-phi}) should give rise to
Floer homology groups $\HF^*(M,\om,\mu,f)$
that are independent of $H$ and $J$.
One might hope that these invariants
can be used to distinguish equivariant 
Hamiltonian isotopy classes. 
In the standard case such results were established
by Seidel~\cite{SEIDEL}. 


\subsection{Boundary value problems}

It would be interesting to consider boundary value
problems for the equations~(\ref{eq:new-ham}) 
or~(\ref{eq:floer}).  The relevant boundary data 
would then involve $\G$-invariant 
Lagrangian submanifolds.  In particular, this should 
lead to Floer homology groups $\HF^*(M,\om,\mu,L_0,L_1)$,
where $L_0$ and $L_1$ are $\G$-invariant Lagrangian
submanifolds of $\mu^{-1}(0)$, 
the critical points are $\G$-orbits 
of intersections of $L_0$ and $L_1$,
and the connecting orbits are solutions of~(\ref{eq:floer})
on the strip $\R\times[0,1]$ that satisfy the 
boundary condition
\begin{equation}\label{eq:bc}
     u(s,0)\in L_0,\qquad 
     u(s,1)\in L_1.
\end{equation}
(See Floer~\cite{FLOER3,FLOER4}, Oh~\cite{OH}, and
Lazzarini~\cite{L} for the standard case.)

\begin{lemma}
Let $L\subset M$ be a connected $\G$-invariant
Lagrangian submanifold.  Then there exists 
a central element $\tau\in\g$ such that
$
     L\subset\mu^{-1}(\tau).
$
\end{lemma}

\begin{proof}
For every $x\in L$ we have
$$
     \im L_x\subset T_xL \subset \ker d\mu(x).
$$
The last inclusion follows from
the fact that $T_xL$ is a Lagrangian subspace
of $T_xM$ and the kernel of $d\mu(x)$ is the symplectic
complement of the image of $L_x$.  The inclusion
$T_xL\subset\ker d\mu(x)$ shows that $\mu$ 
is constant on $L$.  The inclusion 
$\im L_x\subset \ker d\mu(x)$ shows that 
$$
     [\xi,\mu(x)]=d\mu(x)L_x\xi=0
$$
for $\xi\in\g$ and $x\in L$.  
Hence $\mu(x)$ is in the center of $\g$ 
for every $x\in L$.
\end{proof}

If $0$ is a regular value of $\mu$ then the
{\it equivariant diagonal} 
$$
     \Delta^\mu 
     = \left\{(x,gx)\,|\,
     x\in M,\,g\in\G,\,\mu(x)=0\right\}
$$
is a Lagrangian submanifold of 
$\Hat M=M\times M$, with the symplectic
form $\Hat\om=(-\om)\times\om$, and 
is invariant under the action of 
$\Hat\G=\G\times\G$.  The Floer
homology groups of a symplectomorphism
$f:M\to M$ that satisfies~(\ref{eq:mu-phi})
should be isomorphic to the Floer 
homology groups of the Lagrangian pair
$(\Delta^\mu,\Gamma^\mu(f))$ in $\Hat M$,
where $\Gamma^\mu(f)$ is the {\it equivariant
graph} of $f$, i.e. the image 
of $\Delta^\mu$ under $\id\times f$.


\subsection{Adiabatic limits}

That the present theory is, for regular quotients, 
equivalent to the standard theory in $M\dslash\G$
follows from an adiabatic limit argument involving 
the equation
\begin{equation}\label{eq:floer-eps}
     \p_su+J(\p_tu+X_\Psi(u)-X_{H_t}(u)) = 0,\qquad
     \p_s\Psi + \eps^{-2}\mu(u) = 0.
\end{equation}
In the limit $\eps\to 0$ the 
solutions of~(\ref{eq:floer-eps})
degenerate to Floer gradient 
lines in the quotient $M\dslash\G$. 
The details are analogous to the proof
of the Atiyah--Floer conjecture 
in~\cite{DOSAL2,DOSAL3,SAL4}
and to the proof of 
Conjecture~\ref{con:adiabatic} in~\cite{GAIO}.
The resulting theorem should be the existence
of a natural isomorphism
$$
     \HF^*(M,\om,\mu,f)
     \cong \HF^*(M\dslash\G,\bar\om,\bar f),
$$
whenever $0$ is a regular value of $\mu$
and $\G$ acts freely on $\mu^{-1}(0)$. 
Here $\bar\om$ denotes the induced symplectic
form and $\bar f$ the induced
symplectomorphism on $M\dslash\G$. 
In the Lagrangian case there should be a natural 
isomorphism
$$
     \HF^*(M,\om,\mu,L_0,L_1)
     \cong \HF^*(M\dslash\G,\bar\om,\bar L_0,\bar L_1),
$$
where $\bar L_i=L_i/\G\subset M\dslash\G$ for $i=0,1$.


\section{Examples}\label{sec:ex}


\subsection{Vortex equations}\label{sec:vortex}

Consider the standard action of $\G=S^1$ on $M=\C$.
Then a moment map is given by 
\begin{equation}\label{eq:momC}
     \mu(z) = -\frac{i}{2}\left|z\right|^2
\end{equation}
and the quotient space is a point.  Nevertheless,
the space of solutions of~(\ref{eq:new}) is interesting.
Let $(\Sigma,J_\Sigma,\dvol_\Sigma)$ be a compact
Riemann surface, and $P\to\Sigma$ be a circle
bundle of degree $d$.  An equivariant function
$\Theta:P\to\C$ can then be interpreted as a 
section of the line bundle 
$
     E=P\times_{S^1}\C\to\Sigma,
$
a connection $A\in\Aa(P)$ determines a Cauchy-Riemann 
operator 
$$
     \bar\p_A:\Cinf(\Sigma,E)\to\Om^{0,1}(\Sigma,E),
$$
and the equations~(\ref{eq:new}), with
$\mu$ replaced by $\mu+i\tau$ for some $\tau\in\R$,
have the form
\begin{equation}\label{eq:vortex}
     \bar\p_A\Theta=0,\qquad
         *iF_A+\frac{\left|\Theta\right|^2}{2} = \tau.
\end{equation}
These are the {\bf vortex equations}.
The necessary condition~(\ref{eq:positive}) for the 
existence of solutions has the form
$$
     \tau  > \frac{2\pi d}{\Vol(\Sigma)}.
$$
In this case the moduli space is smooth and,
by Proposition~\ref{prop:compact}, it is compact.
These observations are well known~\cite{OGP} as is the fact that
the moduli space 
$$
     \Mm_d(\Sigma)
         = \frac{\left\{(\Theta,A)\,|\,(\ref{eq:vortex})\right\}}
                {\Map(\Sigma,S^1)}
$$
can be identified with the symmetric product 
$
     S^d\Sigma=\Sigma\times\cdots\times\Sigma/S_d
$
(via the zeros of $\Theta$). 
Hence the invariants~(\ref{eq:invariant}) should 
be expressable in terms of the cohomology of $S^d\Sigma$. 
Note that the adiabatic limit argument
of Section~\ref{sec:adiabatic} can, in this case,
be rephrased in the form $\tau\to\infty$ 
(by rescaling $\Theta$) 
and this limit corresponds precisely to the 
argument of Taubes in~\cite{TAUBES2} for the 
Seiberg--Witten equations.


\subsection{Bradlow pairs}\label{sec:bradlow}

Another example with a trivial quotient is the
action of $\G=\U(2)$ on $M=\C^2$.  
A moment map is given by 
$
     \mu(z) = -izz^*/2.
$
Hence the quotient at any nonzero central element
of $\u(2)$ is the empty set. 
Let $P\to\Sigma$ be a principal $\U(2)$-bundle 
of degree 
$
     d = \inner{c_1(E)}{[\Sigma]}
$ 
and consider the Hermitian rank-$2$ bundle
$
     E=P\times_{\U(2)}\C^2 \to\Sigma.
$
Fix a constant 
$
     \tau > \pi d/\Vol(\Sigma)
$
and replace the moment map by $\mu+i\tau\1$. 
Then~(\ref{eq:new}) takes the form
\begin{equation}\label{eq:bradlow}
     \bar\p_A\Theta=0,\qquad
     *iF_A+\frac12\Theta\Theta^* = \tau\1,
\end{equation}
where $A\in\Aa(E)$ and $\Theta\in\Cinf(\Sigma,E)$. 
The moduli spaces 
$$
     \Mm_\tau
         = \frac{\left\{(\Theta,A)\,|\,(\ref{eq:bradlow})\right\}}
                {\Gg(E)}
$$
were studied in detail by Bradlow et al~\cite{BRADLOW,BD,THADDEUS}.  
The invariants~(\ref{eq:invariant})
and the wall crossing numbers should, in this case,
be related to the work of Thaddeus~\cite{THADDEUS}.
He studied the cohomology of the moduli space 
of flat $\U(2)$-connections over $\Sigma$ via Bradlow pairs.
For $\tau$ close to $\pi d/\Vol(\Sigma)$, and large $d$,
$\Mm_\tau$ is a bundle over the moduli space of flat 
$\U(2)$-connections with projective spaces as fibres,
and if $\tau$ is large then $\Mm_\tau$ can be 
identified with a complex projective space. 
The critical parameters are 
$$
      \tau_k = \frac{2\pi k}{\Vol(\Sigma)},\qquad
          \frac{d}{2} < k\le d.
$$
For $\tau=\tau_k$ there are reducible solutions
of~(\ref{eq:bradlow}), i.e. 
$\Mm_\tau$ is not smooth and its singular part 
can be identified with the symmetric product $S^{d-k}\Sigma$.
In the context of this paper it is useful to
recall the following construction~\cite{BRADW}.
Fix a point $z_0\in\Sigma$ and denote by 
$\Gg_0\subset\Map(\Sigma,S^1)$ the codimension-$1$
subgroup of all maps of the form $g=g_0\exp(\xi)$
where $g_0:\Sigma\to S^1$ satisfies $d^*({g_0}^{-1}dg_0)=0$
and $g_0(z_0)=1$ and $\xi:\Sigma\to i\R$ has
mean value zero. Then the quotient space
\begin{equation}\label{eq:Mbradlow}
      \Mm
          = \frac{\left\{(\Theta,A)\,|\,
            \exists\tau > \pi d/\Vol(\Sigma)
            \mbox{ s.t. }(\ref{eq:bradlow})
            \mbox{ holds}\right\}}
            {\left\{g\in\Gg(E)\,|\,\det\circ g\in\Gg_0\right\}}
\end{equation}
is a smooth manifold.  
It carries a Hamiltonian $S^1$-action
with moment map
\begin{equation}\label{eq:mom-bradlow}
      \Mm\to i\R:
      (\Theta,A)\mapsto 
      -\frac{i}{2}\int_\Sigma\left|\Theta\right|^2\dvol_\Sigma.
\end{equation}
Hence $\Mm_\tau$ can be identified with the quotient 
$\Mm\dslash S^1(i(2\pi d-2\tau\Vol(\Sigma))$.


\subsection{Holomorphic curves in projective space}

Consider the standard action of $\G=S^1$ on $\C^{n+1}$. 
Then a moment map is again given by~(\ref{eq:momC}) 
and~(\ref{eq:new}) has the form
\begin{equation}\label{eq:multi-vortex}
     \bar\p_A\Theta_\nu=0,\qquad
         *iF_A+\sum_{\nu=0}^n\frac{\left|\Theta_\nu\right|^2}{2} = \tau,
\end{equation}
where $E\to\Sigma$ is a Hermitian line bundle,
$A\in\Aa(E)$, and $\Theta_0,\dots,\Theta_n\in\Cinf(\Sigma,E)$.
By Proposition~\ref{prop:compact}, the moduli space
of solutions of~(\ref{eq:multi-vortex}) is compact
and transversality can be easily achieved. 
By Conjecture~\ref{con:adiabatic}, the resulting 
invariants~(\ref{eq:invariant})
agree with the Gromov--Witten invariants of $\C P^n$.
However, in contrast to those, they are defined 
in terms of compact smooth moduli spaces. 


\subsection{Toric varieties}

The situation is similar for K\"ahler manifolds that arise
as quotients of $\C^N$ by a subgroup $\G\subset\U(N)$.
Here a moment map is given by 
$$
     \mu(z) = \pi\left(-\frac{i}{2}zz^*\right)
$$
where $\pi:\u(N)\to\g$ denotes the adjoint of the 
inclusion $\g\INTO\u(N)$.  Let $P\to\Sigma$ be
a principal $\G$-bundle and denote by 
$$
     E = P\times_\G\C^N
$$
the associated vector bundle.  
Then the equations~(\ref{eq:new}) 
can be interpreted as equations for a pair
$(\Theta,A)\in\Cinf(\Sigma,E)\times\Aa(E)$
and they have the form
\begin{equation}\label{eq:toric}
     \bar\p_A\Theta = 0,\qquad
         *F_A + \pi\left(-\frac{i}{2}\Theta\Theta^*\right) = \tau
\end{equation}
for some central element $\tau\in\g$.  
Again, Proposition~\ref{prop:compact} guarantees
that the moduli space is compact whenever the moment 
map is proper. 
One gets integer invariants which should correspond
to the Gromov--Witten invariants of the quotient
$
     X = \C^N\dslash\G(\tau).
$
This is interesting, because there are many examples
where $X$ contains holomorphic spheres with negative Chern
number and in these cases the direct definition of the 
Gromov--Witten invariants of $X$ has so far only 
been established over the rationals~\cite{FO,LIT,RUAN}.


\subsection{The Grassmannian and the Verlinde algebra}\label{sec:grass}

In~\cite{WITTEN} Witten conjectured a relation between
the Gromov--Witten invariants of the Grassmannian~\cite{BDW}
and the Verlinde algebra~\cite{VERLINDE,BL}.
For the quantum cohomology ($3$-punctured spheres)
this conjecture was confirmed by Agnihotri~\cite{AGNIHOTRI}.
The Grassmannian can be expressed as a symplectic quotient
$$
     {\rm Gr}(k,n)\cong\C^{k\times n}\dslash\U(k).
$$
Think of $\Theta\in\C^{k\times n}$ as a $k$-frame in $\C^n$.  
If $\Theta$ has rank $k$ then the orthogonal 
complement of its kernel is a $k$-dimensional
subspace of $\C^n$. The group $\U(k)$ acts on
$\C^{k\times n}$ on the left and
the function $\mu:\C^{k\times n}\to\u(k)$ given by
$$
     \mu(\Theta) = -\frac{i}{2}\Theta\Theta^*
$$
is a moment map.  
Thus $\mu^{-1}(-i/2)$ is the space of unitary
$k$-frames in $\C^n$ and its quotient by $\U(k)$
is the Grassmannian.  Now let $P\to\Sigma$
be a principal $\U(k)$-bundle of degree $d$ 
and denote by 
$$
     E=P\times_{\U(k)}\C^k\to\Sigma
$$
the associated complex rank-$k$ bundle.
Fix a real number $\tau>2\pi d/k\Vol(\Sigma)$
and replace $\mu$ by $\mu+i\tau\1$. 
Then~(\ref{eq:new}) takes the form
\begin{equation}\label{eq:grass}
     \bar\p_A\Theta_\nu = 0,\qquad
     *iF_A + \frac12\sum_{\nu=1}^n\Theta_\nu{\Theta_\nu}^*
     = \tau\1,
\end{equation}
where $A\in\Aa(E)$ and $\Theta_1,\dots,\Theta_n\in\Cinf(\Sigma,E)$.
In this case~(\ref{eq:grass}) has reducible solutions
whenever
$$
     \tau = \frac{2\pi d_0}{k_0\Vol(\Sigma)},\qquad
         0<k_0<k,\qquad
         \frac{k_0}{k}d<d_0\le d. 
$$
Thus the moduli space is regular for 
$\tau>2\pi d/\Vol(\Sigma)$ and for these 
values of $\tau$ Conjecture~\ref{con:adiabatic}
asserts that the invariants obtained from the solutions
of~(\ref{eq:grass}) can be identified
with the Gromov--Witten invariants 
of the Grassmanian. On the other hand the opposite 
adiabatic limit $\eps\to\infty$ in~(\ref{eq:rita})
should give rise to an identification with the
invariants of moduli spaces of flat connections 
that appear as the structure constants in the 
Verlinde algebra (cf. Thaddeus~\cite{THADDEUS}). 
Thus the solutions of~(\ref{eq:grass}) might give rise
to a geometric approach for the proof of Witten's 
conjecture. 


\subsection{Anti-self-dual Yang--Mills equations}

There are interesting cases where the solutions of~(\ref{eq:new})
give rise to finite dimensional moduli spaces even though
the symplectic manifold $(M,\om)$ is infinite dimensional.
As an example consider the case of a principal 
bundle $Q\to S$ over a compact oriented Riemann 
surface $S$ with structure group $\SU(2)$ or $\SO(3)$. 
In section~\ref{sec:gauge} we have seen that the space 
$
     M = \Aa(Q)
$
of connections on $Q$ carries a natural symplectic 
structure and that the action of the 
identity component of the gauge group 
$
     \G = \Gg_0(Q)\subset\Gg(Q)
$
is Hamiltonian with moment map 
$\Aa(Q)\to\Lie(\Gg_0(Q)):A\mapsto *F_A$. 
Hence the equations~(\ref{eq:new-loc}),
in local holomorphic coordinates on $\Sigma$, 
have the form
\begin{equation}\label{eq:asd}
\begin{array}{rcl}
     \p_sA-d_A\Phi + *(\p_tA-d_A\Psi) & =& 0,  \\
     \p_s\Psi-\p_t\Phi+[\Phi,\Psi] + \lambda^2*F_A &= & 0,
\end{array}
\end{equation}
where $A(s,t)\in\Aa(Q)$ and 
$\Phi(s,t),\Psi(s,t)\in\Cinf(S,\ad(Q))$
and the metric on $\Sigma$ is $\lambda^2(ds^2+dt^2)$. 
These are the anti-self-dual 
Yang-Mills equations over the product $\Sigma\times S$.
The function $\C\to\Aa(Q):s+it\mapsto A(s,t)$
plays the role of the map $u:\C\to M$ 
in~(\ref{eq:new-loc}).
The symplectic quotient 
$$
     \Mm_Q : = \Aa^\FLAT(Q)/\Gg_0(Q) = M\dslash\G 
$$
is the moduli space of flat connections on $Q$.
It is a symplectic manifold of dimension $6g-6$,
where $g$ is the genus of $S$.  
If $Q$ is an $\SO(3)$-bundle with nonzero second
Stiefel-Whitney class then the moduli space 
$\Mm_Q$ is smooth.  The adiabatic
limit argument of Conjecture~\ref{con:adiabatic} 
here gives rise to a correspondence between
anti-self-dual instantons over $\Sigma\times S$
and holomorphic curves $\Sigma\to\Mm_Q$.  
This is the basic idea of the proof of the 
Atiyah--Floer conjecture~\cite{DOSAL1,DOSAL2,DOSAL3,SAL4}. 
Another reference for this adiabatic limit is the recent 
thesis by Handfield~\cite{HANDFIELD}. 

\begin{remark}
{\bf (i)}
An automorphism $f:Q\to Q$ (that descends to a 
diffeomorphism of $S$) determines an equivariant
symplectomorphism 
$$
     \Aa(Q)\to\Aa(Q):A\mapsto f^*A.
$$
The corresponding isomorphism of the gauge group
is given by 
$$
     \Gg_0(Q)\to\Gg_0(Q):g\mapsto g\circ f.
$$
That the symplectic Floer homology groups of the 
induced symplectomorphism of $\Mm_Q$ are isomorphic
to the instanton Floer homology groups of the corresponding
$3$-dimensional mapping torus was proved in~\cite{DOSAL3}.

\smallskip
\noindent{\bf (ii)}
There are interesting Lagrangian submanifolds of
$\Aa(Q)$ whenever $S$ is the boundary 
of a compact $3$-manifold $Y$ and $Q$ admits a trivialization.
Then the bundle extends over $Y$ and 
the flat connections on $Y$ determine a $\Gg(Q)$-invariant
Lagrangian submanifold of $\Aa(Q)$ (that is contained in
the subset of flat connections).  The general
Atiyah--Floer conjecture~\cite{ATIYAH} relates the 
Floer homology groups of Lagrangian intersections
in $\Mm_Q$, corresponding to two bordisms
$Y_0$ and $Y_1$, to the instanton Floer homology groups
of the closed $3$-manifold $Y=Y_0\cup_SY_1$, whenever the latter
is a homology $3$-sphere. 
\end{remark}


\subsection{Seiberg--Witten equations}

Another infinite dimensional example 
is the space 
$$
     M = 
     \left\{(\Theta,A)\in\Cinf(S,E)\times\Aa(E)\,|\,
     \bar\p_A\Theta = 0
     \right\},
$$
where $E\to S$ is a Hermitian line bundle 
of degree $d=\inner{c_1(E)}{[S]}$ 
over a compact oriented Riemann surface $S$. 
The symplectic form is given by~(\ref{eq:Om})
and Proposition~\ref{prop:Moment} asserts that 
the action of the gauge group 
$
     \G = \Map(S,S^1)
$
on this space is Hamiltonian with moment map
$$
     (\Theta,A)\mapsto *F_A-\frac{i}{2}\left|\Theta\right|^2.
$$
The symplectic quotient $M\dslash\G(-i\tau)$ is the 
moduli space $\Mm_d(S)$ of solutions to the vortex 
equations~(\ref{eq:vortex}) and hence can be 
identified with the $d$-fold symmetric product of $S$.
The equations~(\ref{eq:new-loc}) have 
the form
\begin{equation}\label{eq:SW1}
\begin{array}{rcl}
     \bar\p_A\Theta &= & 0, \\
     \p_s\Theta + \Phi\Theta
     + i(\p_t\Theta+\Psi\Theta) &= & 0, \\
     \p_sA-d\Phi + *(\p_tA-d\Psi) & = & 0,  \\
     \p_s\Psi-\p_t\Phi + 
     \lambda^2(*F_A -i\left|\Theta\right|^2/2 +i\tau)&= & 0,
\end{array}
\end{equation}
where $A(s,t)\in\Aa(E)$, $\Theta(s,t)\in\Cinf(S,E)$,
and $\Phi(s,t),\Psi(s,t)\in\Cinf(S,i\R)$. 
These are the Seiberg--Witten equations over the product
$\Sigma\times S$, so long as the complex structure
on $S$ is independent of $s$ and $t$
(the {\it integrable case}).   More precisely,
the first two equations in~(\ref{eq:SW1}) 
correspond to the Dirac equation and the last two 
to the curvature equation.  The spinor bundle
is a rank-$2$ bundle over $\Sigma\times S$
which naturally splits into a direct sum of two line 
bundles.  In the integrable case one of the 
two components of the spinor vanishes~\cite{WITTEN1} 
and this leads to the simpler form of the 
Seiberg-Witten equations stated above. 

The adiabatic limit argument of Conjecture~\ref{con:adiabatic}
now gives rise to a correspondence between the 
Seiberg--Witten equations over the product
$\Sigma\times S$ and holomorphic curves 
from  $\Sigma$ into the $d$-fold symmetric product
of $S$~\cite{SAL2}.  There is a somewhat more complicated
version of this argument which also applies to the case
where the complex structure on $S$ depends 
on $s$ and $t$.   Then the moduli
spaces of solutions of the vortex equations 
form a bundle over the Teichm\"uller space of $S$,
this bundle carries a natural connection,
and this connection is related in an interesting way
to the full version of the Seiberg--Witten equations
in the nonintegrable case, whenever the $4$-manifold
in question is a fibration with fibre $S$. 
This is discussed in detail in~\cite{SAL2}.
The correspondence between holomorphic curves
and Seiberg--Witten equations indicated here is
different from the one in the work of 
Taubes~\cite{TAUBES1,TAUBES2,TAUBES3,TAUBES4}
where he directly compares the Seiberg--Witten monopoles
over a general symplectic $4$-manifold $X$ 
with holomorphic curves in $X$.  It is likely that 
the two approaches are related via the work of 
Donaldson~\cite{DONALDSON3} on symplectic
Lefschetz fibrations (see also Auroux~\cite{AUROUX}).
Donaldson proved that every symplectic $4$-manifold,
after blowup, admits the structure of a symplectic Lefschetz
fibration 
$$
     X\to S^2
$$
with generic fibre $S$. 
Cutting out the singular fibres 
one obtains a $4$-manifold $W$, fibred over the 
punctured sphere, with cylindrical ends corresponding
to the mapping tori of Dehn twists.  
The adiabatic limit argument of Conjecture~\ref{con:adiabatic}
now relates the Seiberg--Witten monopoles over $X$ 
to holomorphic sections of the bundle $X^{(d)}$,
where the fibres are replaced by the $d$-fold
symmetric products of $S$.  The latter 
correspond to multivalued sections of the bundle
$X\to S^2$.  That these in turn should correspond 
to holomorphic curves in $X$ itself is the 
subject of a current research project by
Donaldson and Ivan Smith. 
The adiabatic limit argument for $W^{(d)}$ 
is the Seiberg--Witten analogue of the
Atiyah--Floer conjecture~\cite{SAL2,SAL4}.
In the $3$-dimensional case this is related 
to the work of Meng--Taubes~\cite{MT},
Hutchings--Lee~\cite{HL1,HL2}, Turaev~\cite{TURAEV},
and Donaldson~\cite{DONALDSON4}.

\begin{remark}
Since there is a correspondence between Donaldson
invariants and holomorphic curves in the moduli space 
$\Mm^\FLAT(S)$ of flat $\SO(3)$-con\-nec\-tions 
over $S$ on the one hand, 
and between the Seiberg--Witten invariants and holomorphic
curves in the symmetric product $\Mm_d(S)$ on the other hand,
it would be interesting to compare the Gromov--Witten
invariants of $\Mm_d(S)$ with those of $\Mm^\FLAT(S)$.
Such a comparison should be related to the picture of 
Thaddeus~\cite{THADDEUS} for the ordinary cohomology
of these spaces, and hence to the study of holomorphic 
curves in the moduli spaces of Bradlow pairs. 
Results in this direction might provide an alternative 
approach (to the one by Pidstrigach--Tyurin~\cite{PT}) 
for the comparison of the Donaldson and the Seiberg--Witten
invariants in the symplectic case.  
The discussion of Section~\ref{sec:bradlow}
shows that this fits into the framework of the 
invariants~(\ref{eq:invariant}).  
To be more precise, equations~(\ref{eq:new-loc})
with target space $\Mm$ given by~(\ref{eq:Mbradlow})
and moment map~(\ref{eq:mom-bradlow})
take the form
\begin{equation}\label{eq:new-bradlow}
\begin{array}{rcl}
    \bar\p_A\Theta &= &0, \\
    \p_s\Theta+\Phi\Theta + i(\p_t\Theta+\Psi\Theta) &= & 0, \\
    \p_sA-d_A\Phi + *(\p_tA-d_A\Psi) &= & 0, \\
    (2\Vol(S))^{-1}\int_S \rm{tr}(\p_s\Psi-\p_t\Phi)\dvol_S
    +\lambda^2(*F_A - i\Theta\Theta^*/2 + i\tau)
    & = & 0,
\end{array}
\end{equation}
where $E\to S$ is a Hermitian rank-$2$ bundle,
$A(s,t)\in\Aa(E)$, $\Theta(s,t)\in\Cinf(S,E)$,
and $\Phi(s,t),\Psi(s,t)\in\Cinf(S,\End(E))$. 
One can prove, with the same techniques as in 
Proposition~\ref{prop:compact}, that the $L^2$-norm
of $\Theta(s,t)$ over $S$ satisfies a universal upper bound
for every solution of~(\ref{eq:new-bradlow})
(over a compact Riemann surface).
Working with~(\ref{eq:new-bradlow}),
instead of holomorphic curves in $\Mm_\tau$ 
(see Section~\ref{sec:bradlow}),
eliminates the problems arising from holomorphic
spheres with negative Chern number,
which exist in $\Mm_\tau$ but not in $\Mm$. 
On the other hand, care must be taken with the 
solutions of~(\ref{eq:new-bradlow}) that satisfy
$\Theta=0$.  
\end{remark}

\smallbreak


\section*{Acknowledgement}

We would like to thank Paul Biran, Simon Donaldson, 
Shaun Martin, and Paul Seidel for enlightening discussions.



\begin{thebibliography}{99}

\small

\bibitem{AGNIHOTRI} S.~Agnihotri,
      Quantum cohomology and the Verlinde algebra,
      PhD thesis, Oxford, 1995.

\bibitem{ATIYAH} M.F.~Atiyah,
      New invariants of three and four dimensional manifolds,
      {\it Proc. Symp. Pure Math.} {\bf 48} (1988).

\bibitem{AB} M.F.~Atiyah and R.~Bott,
      The Yang--Mills equations over Riemann surfaces,
      {\it Phil. Trans. R. Soc. Lond. A} {\bf 308} (1982), 523--615.

\bibitem{AUROUX} D.~Auroux,
      Asymptotically holomorphic families of symplectic submanifolds,
      {\it Geom. Funct. Anal.} {\bf 7} (1997), 971--995.

\bibitem{BDW} A.~Bertram, G.~Daskalopoulos, and R.~Wentworth,
      Gromov invariants for holomorphic maps from 
      Riemann surfaces to Grassmannians,  Preprint, April 1993.

\bibitem{BL} J.M.~Bismut and F.~Labourie,
      Symplectic geometry and the Verlinde formula,
      Preprint, Universit\'e de Paris Sud, Orsay, 1999. 

\bibitem{BRADLOW} S.~Bradlow,
      Special metrics and stability for holomorphic bundles with 
      global sections, 
      {\it J. Diff. Geom.} {\bf 33} (1991), 169--213.

\bibitem{BD} S.~Bradlow and G.~Daskalopoulos,
      Moduli of stable pairs for holomorphic bundles
      over Riemann surfaces, 
      {\it Int. J. Math.} {\bf 2} (1991), 477--513.

\bibitem{BRADW} S.~Bradlow, G.~Daskalopoulos, and R.~Wentworth,
      Birational equivalence of vortex moduli spaces,
      {\it Topology} {\bf 35} (1996), 731--748.

\bibitem{DAVIES} T.~Davies,
      The Yang-Mills functional over Riemann surfaces
      and the loop group, PhD thesis, Warwick 1996.

\bibitem{DONALDSON1} S.K.~Donaldson,
      Instantons and geometric invariant theory, 
      {\it Comm. Math. Phys.} {\bf 93} (1984), 453--460.

\bibitem{DONALDSON2} S.K.~Donaldson,
      Polynomial invariants of smooth four-manifolds, 
      {\it Topology} {\bf 29} (1990), 257--315.

\bibitem{DONALDSON3} S.K.~Donaldson,
      Symplectic submanifolds and almost complex geometry,
      {\it J. Diff. Geom.} {\bf 44} (1996), 666-705.

\bibitem{DONALDSON4} S.K.~Donaldson,
      Topological field theories and the formulae of 
      Casson and Meng-Taubes, Preprint, 1999.

\bibitem{DONALDSON5} S.K.~Donaldson,
      Symmetric spaces, K\"ahler geometry, and Hamiltonian dynamics,
      Preprint, 1999.

\bibitem{DONALDSON6} S.K.~Donaldson,
      Moment maps and diffeomorphisms,
      Preprint, 1999.

\bibitem{DOSAL1}  S.~Dostoglou and D.A.~Salamon
      Instanton homology and symplectic fixed points,
      in {\it Symplectic Geometry}, edited by D.~Salamon,
      Proceedings of a Conference, 
      LMS Lecture Notes Series {\bf 192},
      Cambridge University Press, 1993, pp. 57--94.

\bibitem{DOSAL2}  S.~Dostoglou and D.A.~Salamon
      Cauchy-Riemann operators, self-duality, and the spectral flow, 
      in {\it First European Congress of Mathematics,
      Volume I, Invited Lectures (Part 1)},
      edited by A.~Joseph, F.~Mignot, F.~Murat,
      B.~Prum, R.~Rentschler,  Birkh\"auser Verlag,
      Progress in Mathematics, {\bf Vol. 119},
      1994, pp. 511--545.

\bibitem{DOSAL3}  S.~Dostoglou and D.A.~Salamon,
      Self-dual instantons and holomorphic curves,
      {\it Annals of Mathematics} {\bf 139} (1994), 581--640.

\bibitem{EGH} Y.~Eliashberg, A.~Givental, and H.~Hofer,
      Contact homology,
      in preparation.

\bibitem{FLOER3} A.~Floer,
      The unregularized gradient flow of the symplectic action,
      {\it Comm. Pure Appl. Math.} {\bf 41} (1988), 775--813.

\bibitem{FLOER4} A.~Floer,
      Morse theory for Lagrangian intersections,
      {\it J. Diff. Geom.} {\bf 28} (1988), 513--547.

\bibitem{FLOER1} A.~Floer,
      An instanton invariant for $3$-manifolds,
      {\it Commun. Math. Phys.} {\bf 118} (1988), 215--240.

\bibitem{FLOER2} A.~Floer,
      Symplectic fixed points and holomorphic spheres,
      {\it Comm. Math. Phys.} {\bf 120} (1989), 575--611.

\bibitem{FO} K.~Fukaya and K.~Ono,
      Arnold conjecture and Gromov--Witten invariants for 
      general symplectic manifolds,
      Preprint, February 1996.

\bibitem{FUKAYA} K.~Fukaya,
      Lectures in Cortona and Kyoto,
      June/July 1999.

\bibitem{GAIO} A.R.~Gaio,
      J-holomorphic curves and moment maps,
      PhD thesis, University of Warwick, 1999.

\bibitem{OGP} O.~Garc\'ia-Prada,
      A direct existence proof for the vortex equations
      over a compact Riemann surface,
      {\it Bull. London Math. Soc.} {\bf 26} (1994), 88--96.

\bibitem{GIVENTAL1} A.B.~Givental,
      Equivariant Gromov--Witten invariants, 
      Pre\-print, Ber\-keley, 1996.

\bibitem{GIVENTAL2} A.B.~Givental,
      Stationary phase integrals, quantum Toda lattices, 
      flag manifolds, and the mirror conjecture,
      Preprint, Berkeley, 1996.

\bibitem{GIVENTAL3} A.B.~Givental,
      A tutorial on quantum cohomology,
      Lecture Notes for the IAS/PCMI Graduate Summer School
      on Symplectic Geometry and Topology,
      December 1997.

\bibitem{GRO} M.~Gromov,
      Pseudo holomorphic curves in symplectic manifolds,
      {\it Invent. Math.} {\bf 82} (1985), 307--347.

\bibitem{GS1} V.~Guillemin and S.~Sternberg,
      Geometric quantization and multiplicities of group 
      representations,
      {\it Invent. Math.} {\bf 67} (1982), 515--538.

\bibitem{GS} V.~Guillemin and S.~Sternberg,
      Birational equivalence in the symplectic category,
      {\it Invent. Math.} {\bf 97} (1998), 485--522.

\bibitem{HANDFIELD} F.G.~Handfield,
      Adiabatic limits of the anti-seld-dual equation,
      PhD thesis, Austin, Texas, 1998.


\bibitem{HOSA}
      H.~Hofer and D.A.~Salamon,
      Floer homology and Novikov rings,
      {\it The Floer Memorial Volume}, 
      edited by H.~Hofer, C.~Taubes, A.~Weinstein, and
      E.~Zehnder, Birkh\"auser 1995, pp 483--524.

\bibitem{HL1} M.~Hutchings and Y.-J.~Lee,
      Circle valued Morse theory, Reidemeister torsion,
      and Seiberg--Witten invariants of $3$-manifolds,
      to appear in {\it Topology}.

\bibitem{HL2} M.~Hutchings and Y.-J.~Lee,
      Circle valued Morse theory and Reidemeister torsion,
      to appear in {\it Math. Research Letters}.

\bibitem{KIRWAN} F.~Kirwan,
      {\it Cohomology of Quotients in Symplectic and 
      Algebraic Geometry}, Princeton University Press, 1984.

\bibitem{KONTSEVICH}  M.~Kontsevich,
      Enumaration of rational curves via torus actions,
      Preprint, 1994.

\bibitem{KOMA} M.~Kontsevich and Yu.~Manin,
      Gromov--Witten classes, quantum cohomology,
      and enumerative geometry,
      Preprint, 1993.

\bibitem{KM} P.~Kronheimer and T.S.~Mrowka,
      The genus of embedded surfaces in the projective plane,
      {\it Math. Res. Letters} {\bf 1} (1994), 797--808.

\bibitem{L} L.~Lazzarini,
      Existence of somewhere injective pseudoholomorphic discs,
      Preprint, 1999. 
          
\bibitem{LT} G.~Liu and G.~Tian,
      Floer Homology and Arnold Conjecture,
      Preprint, August 1996, revised May 1997.

\bibitem{LIT}  Jun Li and Gang Tian,
      Virtual moduli cycles and GW invariant, Preprint 1996.

\bibitem{MARTIN} S.K.~Martin,
      Symplectic geometry and gauge theory,
      PhD thesis, Oxford, 1997. 

\bibitem{MARTIN1} S.K.~Martin,
      Transversality theory, cobordisms,
      and invariants of symplectic quotients, 
      Preprint, IAS, February 1999.

\bibitem{MARTIN2} S.K.~Martin,
      Symplectic quotients by a nonabelian group and 
      by its maximal torus, Preprint, IAS, March 1999.

\bibitem{MARTIN3} S.K.~Martin,
      private communication.

\bibitem{MS1} 
      D.~McDuff and D.~Salamon, 
      {\it $J$-holomorphic Curves and Quantum Cohomology},
      University Lecture Series {\bf 6}, 
      American Mathematical Society, Providence, RI, 1994.  

\bibitem{MS2} D.~McDuff and D.~Salamon,
      {\it Introduction To Symplectic Topology},
      Oxford University Press, 1995, 2nd edition 1998.

\bibitem{MT} G.~Meng and C.Taubes,
      $\SW$ = Milnor Torsion,
      {\it Math. Research Letters} {\bf 3} (1996), 661.

\bibitem{MILNOR}  J.~Milnor,
      {\it Topology From The Differentiable Viewpoint\/},
      University Press of Virginia, Charlottesville, 1965. 
          
\bibitem{MFK} D.~Mumford, J.~Fogarty, and F.~Kirwan,
      {\it Geometric Invariant Theory}, Springer, New York, 1994.

\bibitem{MUNDET} I.~Mundet,
      Yang-Mills-Higgs theory for symplectic fibrations,
      PhD thesis, Madrid, April 1999. 
          
\bibitem{NS} M.S.~Narasimhan and C.S.~Seshadri,
      Stable and unitary vector bundles over compact 
      Riemann surfaces, {\it Ann. Math.} {\bf 82} (1965), 540--564.

\bibitem{OH}  Y.-G.~Oh, 
      Floer cohomology of Lagrangian intersections 
      and pseudoholomorphic discs,
      {\it Comm. Pure Appl. Math.} {\bf 46} (1993), 949--994.

\bibitem{PT} V.~Pidstrigach and A.~Tyurin,
      Localization of Donaldson polynomials along 
      Seiberg--Witten classes,
      Preprint No.~75, Universit\"at Bielefeld, 1995.

\bibitem{PSS} S.~Piunikhin, D.~Salamon, and M.~Schwarz,
      Symplectic Floer-Donaldson theory and quantum cohomology,
      in {\it Contact and Symplectic Geometry}, 
      edited by C.B.~Thomas,
      Publications of the Newton Institute,
      Cambridge University Press 1996, 171--200.

\bibitem{PS} A.~Pressley and G.~Segal,
      {\it Loop Groups}, Oxford University Press, 1986.

\bibitem{RUAN2} Y.~Ruan,
      Topological sigma model and Donaldson type 
      invariants in Gromov theory. 
      {\it Duke Mathematical Journal} 
      {\bf 83} (1996), 461--500

\bibitem{RUAN} Y.~Ruan,
      Virtual neighbourhoods and monopole equations,
      Preprint, March 1996.

\bibitem{RUAN1} Y.~Ruan,
      Surgery, quantum cohomology, and 
      birational geometry, Preprint, 1999.

\bibitem{SAL3} D.A.~Salamon,
      Lagrangian intersections, $3$-manifolds with boundary,
      and the Atiyah--Floer conjecture,
      in {\it Proceedings of the ICM}, Z\"urich,
      1994, Birkh\"auser, Basel, 1995, Vol. 1, 526--536.

\bibitem{SAL2} D.A.~Salamon,
      Seiberg--Witten invariants of mapping tori,
      symplectic fixed points, and Lefschetz numbers,
      {\it Turkish J. of Math.} {\bf 23} (1999), 117--143.

\bibitem{SAL4} D.A.~Salamon,
      Quantum products for mapping tori and the Atiyah--Floer 
      conjecture,  Preprint, ETH-Z\"urich, July 1999.

\bibitem{SAL1} D.A.~Salamon,
      {\it Spin geometry and Seiberg--Witten invariants},
      to appear in Birkh\"auser Verlag.

\bibitem{SEIDEL} P.~Seidel,
      Floer homology and the symplectic isotopy problem,
      PhD thesis, Oxford, 1997.

\bibitem{SEIDEL1} P.~Seidel,
      $\pi_1$  of symplectic automorphism groups 
      and invertibles in quantum cohomology rings.
      {\it Geometric and Functional Analysis}
      {\bf 7} (1997), 1046--95.

\bibitem{TAUBES1} C.H.~Taubes,
      The Seiberg--Witten and the Gromov invariants,
      {\it Math. Res. Letters} {\bf 2} (1995), 221--238.

\bibitem{TAUBES2} C.H.~Taubes,
      $\SW\IMP\Gr$: From the Seiberg--Witten equations 
      to pseudoholomorphic curves,
      {\it J. Amer. Math. Soc.} {\bf 9} (1996), 845--918.

\bibitem{TAUBES3} C.H.~Taubes,
      Counting pseudo-holomorphic submanifolds in dimension $4$,
      Preprint, Harvard, 1996.

\bibitem{TAUBES4} C.H.~Taubes,
      $\Gr\IMP\SW$: From pseudoholomorphic curves
      to the Seiberg--Witten invariants,
      Preprint, Harvard, 1996.

\bibitem{THADDEUS} M.~Thaddeus,
      Stable pairs, linear systems, and the Verlinde algebra,
      {\it Inv. Math.} {\bf 117} (1994), 317--353. 

\bibitem{TURAEV} V.~Turaev,
      Torsion invariants of spin$^c$ structures
      on $3$-manifolds, {\it Math. Res. Letters}
      {\bf 4} (1997), 679--695.

\bibitem{UHL} K.~Uhlenbeck,
      Connections with $L^p$ bounds on the curvature,
      {\it Commun. Math. Phys.} {\bf 83} (1982), 31--42.

\bibitem{VERLINDE} E.~Verlinde,
      Fusion rules and modular transformations
      in 2-D conformal field theory,
      {\it Nucl. Phys. B} {\bf 300} (1988), 360.

\bibitem{VITERBO} C.~Viterbo,
      Functors and computations in Floer homology
      with Applications~I,
      Preprint, 1998.  To appear in {\it GAFA}.

\bibitem{WITTEN} E.~Witten,
      The Verlinde algebra and the quantum cohomology
      of the Grassmannian, Preprint, 
      iassns-hep-93/41, December 1993.

\bibitem{WITTEN1} E.~Witten,
      Monopoles and $4$-manifolds, 
      {\it Math. Res. Letters} {\bf 1} (1994), 769--796.

\end{thebibliography}
\end{document}